# Axial Behaviour of Pre-Damaged RC Short Columns Retrofitted with Square Corrugated Steel Jackets

Ligui Yang[a], Fan Yang[b], Yuyin Wang [c,d], Shuangshuang Jin[a,*], Yong Fang [c,d], Yunwen Chen[a], Qilong Xia[b]

[a] School of Civil Engineering, Chongqing Jiaotong University, Chongqing 400074, China;

[b] Chongqing Construction Engineering Group Co., Ltd;

[c] Key Lab of Structures Dynamic Behaviour and Control of the Ministry of Education, Harbin Institute of Technology, Harbin150090, China;

[d] School of Civil Engineering, Harbin Institute of Technology, Heilongjiang, Harbin150090, China

## Abstract

Existing concrete columns may experience degradation in mechanical performance due to long-term environmental erosion, accidental impacts, or changes in functional requirements. This study proposes a strengthening method employing square corrugated steel jackets as external confinement, which significantly enhances both the bearing capacity and ductility of existing reinforced concrete (RC) columns. Axial compression tests were conducted on ten short column specimens to evaluate the effects of corrugated steel thickness (1.6, 2.0, and 2.7 mm), preloading level before jacketing (40%, 60%, and 100% of the original capacity), and connection type (welding vs. bolting). A computational model was developed to predict the ultimate bearing capacity of the strengthened sections. The main findings are as follows: (1) The corrugated steel jackets increased the ultimate bearing capacity of the existing RC columns by 34.6% to 67.3%. (2) Given the relatively low confinement efficiency in square sections, thinner corrugated steel plates can be used in a material-efficient manner to achieve comparable strengthening effects. (3) Fully welded connections between corrugated plates induce less stress concentration and provide better transverse confinement effectiveness compared to bolted connections. (4) In a pre-unloaded column, greater existed damage causes concrete softening and increased lateral expansion under re-compression. This dilation promotes a tighter interaction between the core concrete and the outer jacket, activating stronger passive confinement after being jacketed. (5) The low longitudinal stresses in the jacket indicate that its primary role is to provide lateral confinement rather than to resist axial loads directly. (6) It is recommended to employ a calculation method that accounts for both pre-damage and confinement effects to ensure a conservative and reliable design of corrugated steel-jacketed RC columns with pre-damage.

| **Nomenclature** | |
|---|---|
| $A_L$ | Area of angle steel ($mm^2$) |
| $A_a$ | Cross sectional areas of equivalent steel tube ($A_a=A_L+A_{cs}$) ($mm^2$) |
| $A_c$ | Concrete cross-sectional area ($mm^2$) |
| $A_{cc}$ | Confined concrete area ($mm^2$) |
| $A_{cs}$ | Areas of corrugated steel ($mm^2$) |
| $A_g$ | Areas of grout ($mm^2$) |
| $A_s$ | Areas of reinforcing bars ($mm^2$) |
| $A_1$ | Total concrete area of confined zone ($mm^2$) |
| $A_2$ | Total concrete area of unconfined zone ($mm^2$) |
| $A_3$ | Total grout area of confined zone ($mm^2$) |
| $A_4$ | Total grout area of unconfined zone ($mm^2$) |
| B | Width of jacket (mm) |
| CFCST | Concrete-filled corrugated steel tube |
| CFDST | Circular-section galvanized corrugated steel tube infilled with concrete |
| CFRP&GFRP | Fiber reinforced plastic |
| CS | Specimen confined by corrugated steel jacket |
| DI | Ductility index |
| E | Elastic modulus (MPa) |
| ECC | Engineered cementitious composites |
| $EI_{eff}$ | The effective flexural stiffness (MPa) |
| $E_a$ | Elastic modular of the structural steel section (MPa) |
| $E'_c$ | Elastic modular of the pre-damaged concrete section (MPa) |
| $E_s$ | Elastic modulus of steel (MPa) |
| $E_s^t$ | Tangential modular (MPa) |
| $I_a$ | Inertia moment of the structural steel section ($mm^4$) |
| $I_c$ | Inertia moments of area of the cracked concrete section ($mm^4$) |
| $I_s$ | Inertia moments of area of the reinforcement for the bending plane being considered ($mm^4$) |
| L | Length of column (mm) |
| N | Load (kN) |
| $N_{cr}$ | Elastic critical force (kN) |
| $N_u$ | Ultimate load of specimens (kN) |
| $N_{u,RC}$ | Ultimate axial capacity of the reinforced concrete specimen/load bearing capacity of RC section (kN) |
| $N_0$ | Preloading/pre-unloading value in the RC section (kN) |
| RC | Reinforced concrete |
| SI | Strength index |
| $S_1$ | Unconfined zone area ($mm^2$) |
| $S_{1,RC}$ | Concrete area in unconfined zone ($mm^2$) |
| UHPC | Ultra-high-performance concrete |
| $b$ | Width of RC section (mm) |
| $d$ | Diameter of steel tube (mm) |
| $f_b$ | Yield strength of confining steel (MPa) |
| $f_c$ | Compressive concrete strength (MPa) |
| $f_c$’ | Residual compressive strength of concrete with pre-damage (MPa) |
| $f_{cc}$ | Confined concrete compressive strength (MPa) |
| $f_{cc}$’ | Compressive strength of pre-damaged and jacketed concrete (MPa) |
| $f_{ck}$ | Characteristic value of the cylinder compressive strength of concrete at 28 days (MPa) |
| $f_{c,g}$ | Compressive strength of grouting material (MPa) |
| $f_{cc,g}$ | Compressive strength of confined grouting material (MPa) |
| $f_{co}$ | Compressive strength of plain concrete (MPa) |
| $f_r$ | Average onfining stress of jacket (MPa) |
| $f_{sc}$ | Compressive strength of the composite section (MPa) |
| $f_y$ | Yield strength (MPa) |
| $f_{y,L}$ | Yield strength of angle steel (MPa) |
| $f_{y,cs}$ | Yield strength of corrugated steel plate (MPa) |
| $f_{y,s}$ | Yield strength of longitudial bars (MPa) |
| $f_u$ | Ultimate strength (MPa) |
| $h_c$ | Corrugation depth (mm) |
| $h_{cs}$ | Corrugation height (mm) |
| $k_d$ | Strength degradation factor of concrete |
| $l_{cs}$ | Corrugation profile length (mm) |
| $r$ | Imaginary radius of unconfined zone (mm) |
| $t$ | Thickness of steel tube (mm) |
| $t_{cs}$ | Thickness of corrugated steel plates (mm) |
| $t_g$ | Thickness of grouting (mm) |
| $\alpha$ | Steel ratio |
| $\Delta$ | Displacement (mm) |
| $\varepsilon_c$ | Compressive strain in the concrete at the peak stress |
| $\varepsilon_h$ | Horizontal strain |
| $\varepsilon_y$ | Yield compressive strain |
| $\varepsilon_v$ | Vertical strain |
| $\varepsilon_0$ | Compressive strain at the start of unloading |
| $\varepsilon_{85\%}$ | Axial compressive strain corresponding to $0.85N_u$ in the descending stage |
| $\xi$ | Confinement factor |
| $\eta_{cs}$ | Amplification coefficient |
| $\mu_s$ | Poisson ratio |
| $\mu_{sp}$ | Elastic-plastic Poisson ratio |
| $\sigma_h$ | Horizontal stress (MPa) |
| $\sigma_v$ | Vertical stress (MPa) |

## 1. Introduction

The safety and durability of buildings are critical concerns in modern urbanization [1]. Structures are vulnerable to damage or collapse due to natural hazards such as wind, floods, fires, and earthquakes (**Fig. 1**), while aging buildings further exhibit progressive deterioration in load-carrying capacity and ductility[2]. Rather than demolishing such structures, repair and reinforcement present a more sustainable and economically viable alternative. This approach is particularly relevant for aging buildings from the early 20th century requiring retrofitting [3],[4].

Columns are the most critical vertical load-bearing elements in building structures, and their load-carrying capacity is essential for ensuring overall structural safety[5]. In high-rise buildings, lower-story columns are subjected to relatively high axial loads and are more susceptible to damage if their bearing capacity is inadequate [6,7], which may lead to significant safety hazards. To address these issues, various reinforcement methods for reinforced concrete columns have been proposed, as illustrated in **Fig. 2**. These include: section enlargement with concrete [8-10], jacketing with ferrocement composites [11], strengthening with steel angles and strips [12,13], and jacketing with steel sheets [14-16].

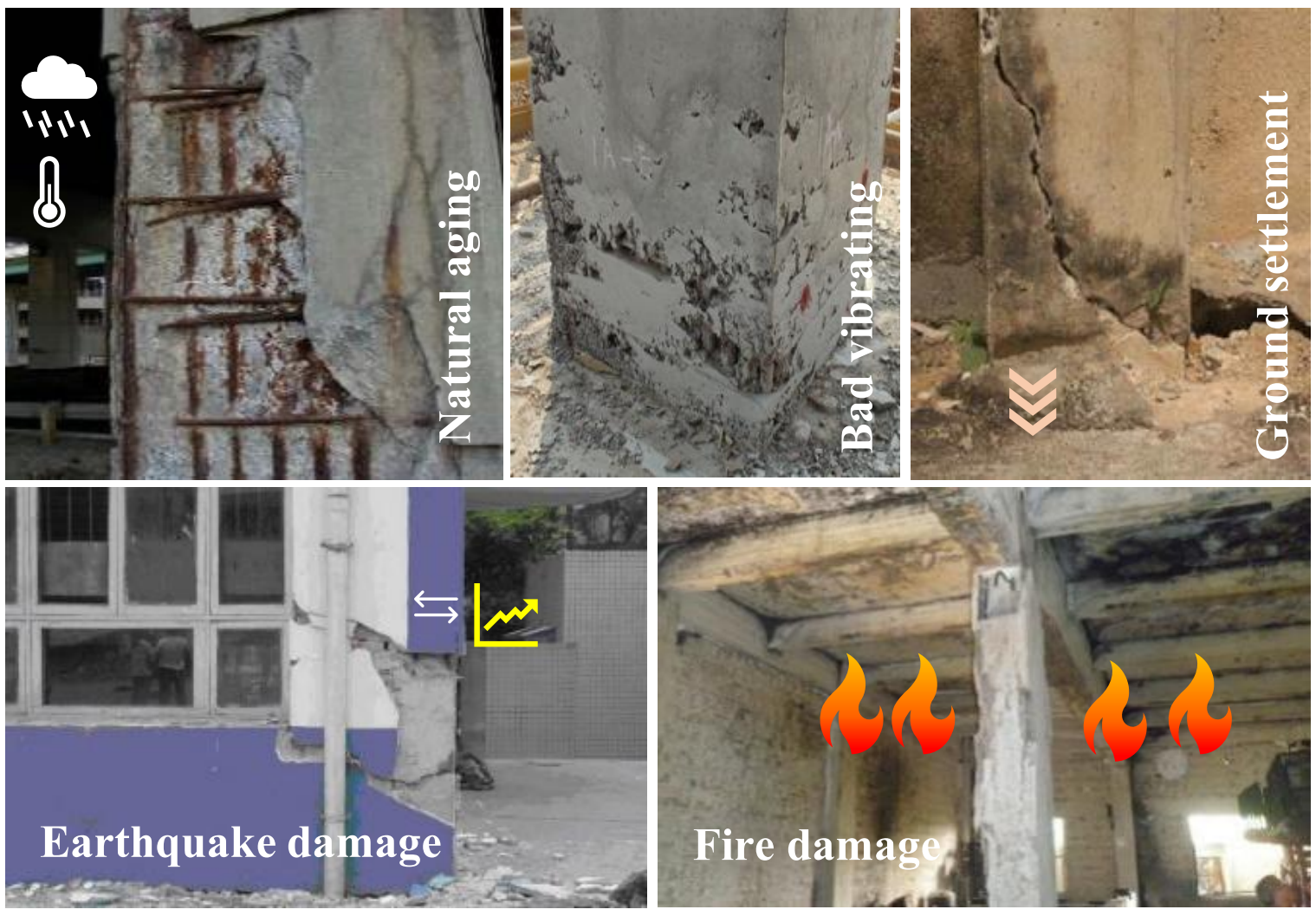

**Fig. 1. Damage and defects of column.**

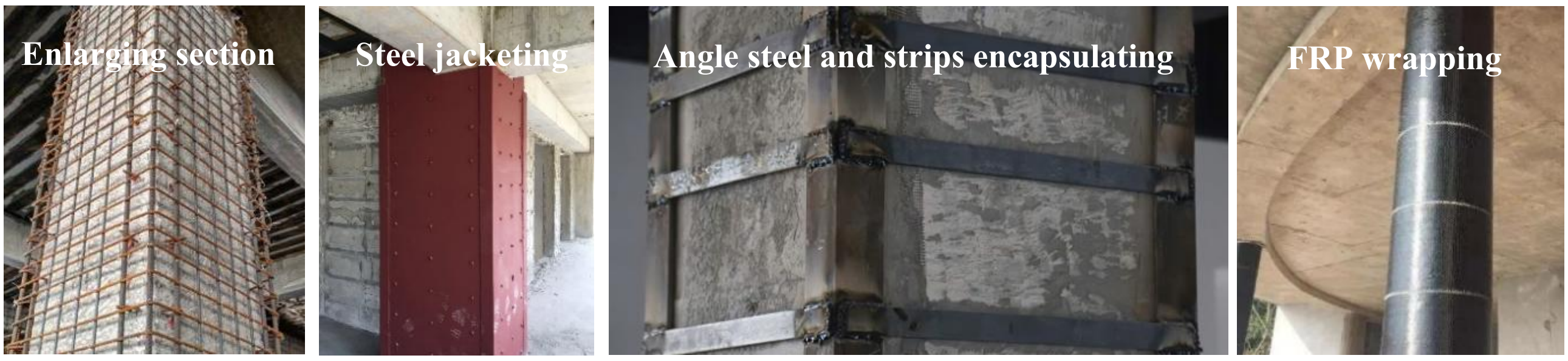

**Fig. 2. Typical strengthening method of column.**

In recent years, more methods have been developed to strengthen reinforced concrete columns using novel high-strength materials [17,18], such as fiber reinforced plastic (GFRP & CFRP) [19-21], ultra-high-performance concrete (UHPC) [22-24], and engineered cementitious composites (ECC)[25-28]. Ghobarah and Bidda et al. (1997) [29,30] proposed a method of reinforcing beams and columns using corrugated steel sleeves. As shown in **Fig. 3**, the corrugated steel sheets are assembled into sleeves to encase the damaged column ends and beam-column joints after an earthquake, which is similar to Tommi's earlier proposal for corrugated steel tube columns [31]. The externally applied corrugated steel jackets significantly increase the load-bearing capacity of the columns due to their circumferential confinement to the core concrete. Yang et al. (2020) [32] investigated the compressive behaviour of the concrete-filled corrugated steel tubular columns, and found that even with thin walls, corrugated steel can provide excellent circumferential confinement to the core component.

Considering the mechanical advantages, concrete columns strengthened with a thin-walled corrugated steel jacket may exhibit good mechanical performance due to the combined benefits of the jacket's encasement and the corrugations' interlocking effect [33]. The jacket also prevents secondary cracking or spalling of the repair layer during service, thereby eliminating the need for frequent inspections. In addition, retrofitting existing reinforced concrete columns with galvanized corrugated steel provides excellent corrosion resistance due to the protective zinc coating and dense oxide film on the surface of the galvanized steel [34]. This results in both maintenance-free service performance and reduced management/maintenance costs. The corrugated steel jacket can serve as formwork due to its significant transverse stiffness. Compared to conventional steel plate reinforcement, this method reduces structural steel usage and minimizes on-site welding operations.

A composite member can be obtained by using corrugated steel plates or tubes to enhance a reinforced concrete (RC) column. For similar concrete columns, Fang et. al. [35,36] investigated the load-bearing capacity of circular-section galvanized corrugated steel tube infilled with concrete (CFDST) and built a design approach for such a novel composite column. For concrete-filled square corrugated steel tubular columns, Zou Y. [37], Yu C. [38], and Nassirnia M. [39] et al. conducted compression and eccentric compression tests, demonstrating that the sleeve method effectively constrains the core concrete. They explained, from a theoretical perspective, why corrugated steel sleeves can significantly restore and enhance the mechanical performance of concrete. Lu and Fang et al. [40] studied the hysteresis performance of concrete-filled double-skin corrugated steel tubes under combined compression and cyclic lateral loads, revealing that the confinement effect of the corrugated steel tube plays a critical role in improving the ductility of the column. They also proposed design formulas for the compression-bending combined resistance of such a hollow-section member. Handousa A. and Abdellatief M. et al. [41] employed two machine learning models to predict the ultimate load-bearing capacity of concrete-filled corrugated steel tubular columns: Artificial Neural Network (ANN) and Gaussian Process Regression (GPR).

As the structural enhancement mechanisms of corrugated steel to the damaged concrete is overlooked in existing studies, this investigation examines the corrugated steel jacketing the concrete column with damage through comprehensive experimental and analytical approaches. This study mainly focuses on the reinforcement of the structure with low axial pre-compression conditions in low-rise buildings, or the retrofitting of unloaded structures in multi-story buildings. Ten short-column specimens were fabricated and subjected to axial compression testing, with particular emphasis on elucidating failure patterns, and ultimate load characteristics. A parametric test was systematically conducted, incorporating three primary variables: connection configuration, corrugated steel plate thickness (ranging 1.6-2.7 mm), and initial damage severity (quantified by 40-100% preloading level). Throughout experimental procedures, the strain evolution in the corrugated steel jackets was continuously monitored, and complemented by elastoplastic analysis to characterize stress redistribution phenomena. The study clarifies the composite action mechanism whereby the confinement from corrugated steel effectively enhances the residual capacity of the damaged concrete core through triaxial compression. Building upon the established ultimate state model of axially compressed sections, a novel theoretical method was developed incorporating confinement and damage reduction factors, enabling precise prediction of retrofitted RC columns' load-bearing capacity with less than 1.86% deviation from experimental values. Calculation methods recommended by the existing codes were modified for such a strengthened cross section. From the systematically comparison results, a simple calculation formula was suggested in the application.

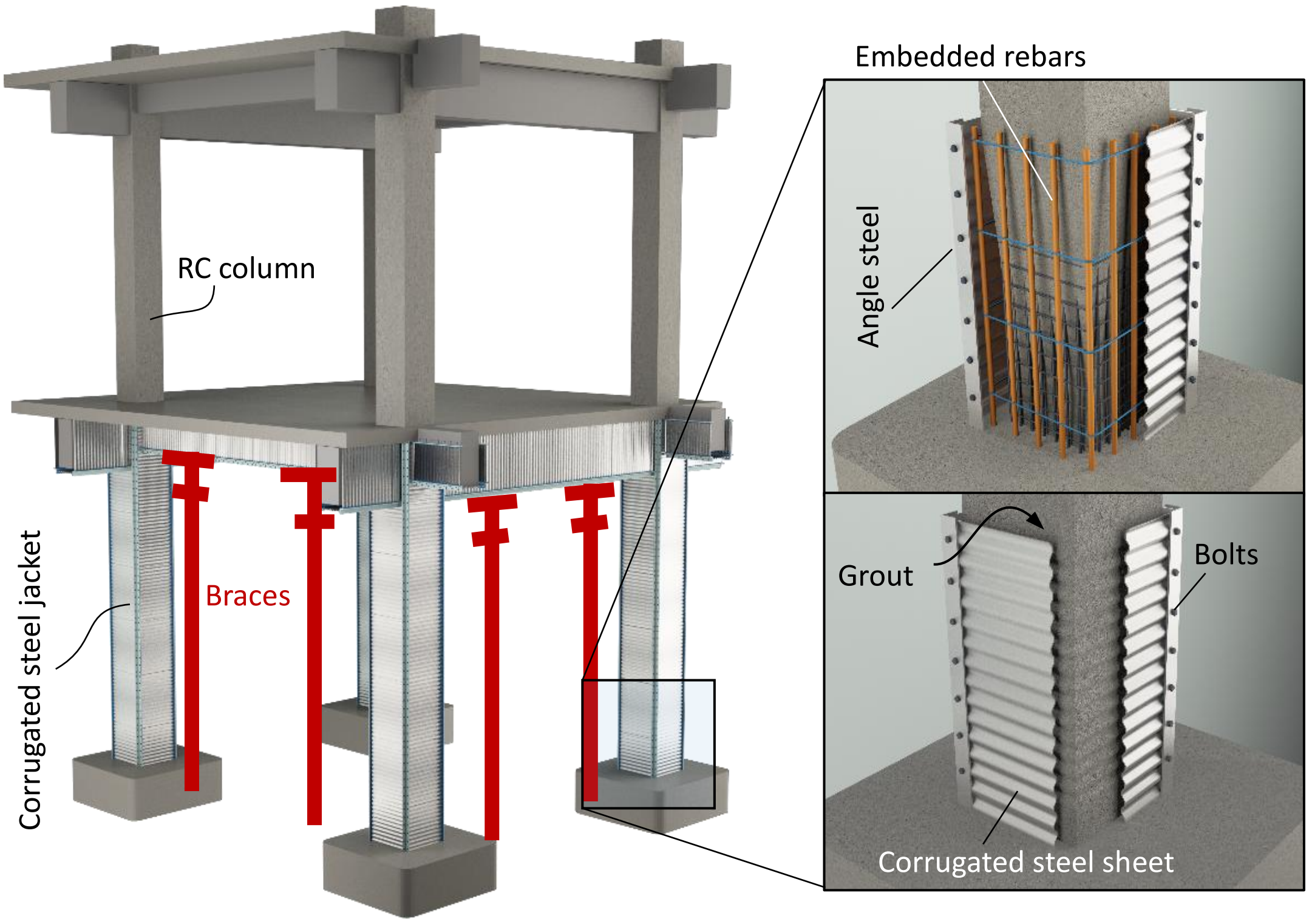

**Fig. 3. Conception of strengthening RC columns by corrugated steel jacket.**

## 2. Experimental arrangements

### 2.1 Specimens parameters

As detailed in **Table 1**, a total of ten specimens were fabricated for this experimental program. The confinement effect of the corrugated steel jackets was primarily evaluated by comparing the performance of specimens with different jacket thicknesses (1.6, 2.0 and 2.7 mm). Considering the potential damage in existing structures, a preload was applied to the reinforced concrete (RC) section to induce plastic damage in the concrete, followed by immediate unloading. This process simulates the bracing effect in actual strengthening operations. As the study focuses on the confinement effect of corrugated steel jackets on pre-damaged concrete, the preloading was conducted as a short-term and quasi-static monotonic cycle, excluding complex concrete creep effects. Preload levels from 40% to 100% of the ultimate axial capacity of the reinforced concrete specimen ($N_{u,RC}$) were adopted to cover a comprehensive range of damage degrees. Nevertheless, strengthening under sustained load, which occurs without bracing or with incomplete unloading, remains a topic for future investigation. For quantitative comparison, an additional plain RC short column was tested to failure to establish a baseline load-displacement curve.

The RC parts of all the specimens are fabricated with the same reinforcement ratio (longitudinal bars ratio=1.1%, and stirrups ratio=0.5%). The corrugated steel plates are formed by bending Q355 steel plates, and the corrugation profile length×height ($l_{cs}$×$h_{cs}$) is 68×13 mm×mm (wave heigh-length ratio = 0.19), which is generally adopt in practice. The gap between the concrete surface and jacket is filled with grout.

**Table 1. Parameters of specimens.**

| Specimens | $t_{cs}$ (mm) | $t_g$ (mm) | Dimension (mm) | $L$ (mm) | $N_0$ | Preload ratio $N_0/N_{u,RC}$ | $f_c$ (MPa) | $N_u$ (kN) | Connection type |
|---|---|---|---|---|---|---|---|---|---|
| RC | / | / | 200×200 | 748 | / | / | 49.7 | 2117.7 | / |
| CS-A-2-0.4 | 2.0 | 14.0 | 241×241 | 748 | 880 | 0.4 | 49.7 | 2957.4 | Connection A |
| CS-B-2-0.4 | 2.0 | 21.5 | 256×256 | 748 | 880 | 0.4 | 49.7 | 3101.1 | Connection B |
| CS-C-2-0.4 | 2.0 | 17.0 | 247×247 | 748 | 880 | 0.4 | 49.7 | 3125.6 | Connection C |
| CS-D-1.6-0.4 | 1.6 | 22.0 | 257×257 | 748 | 880 | 0.4 | 49.7 | 3361.3 | Connection D |
| CS-D-2.7-0.4 | 2.7 | 22.0 | 257×257 | 748 | 880 | 0.4 | 49.7 | 3314.0 | Connection D |
| CS-D-2-0.4 | 2.0 | 22.0 | 257×257 | 748 | 880 | 0.4 | 49.7 | 3543.4 | Connection D |
| CS-D-2-0.6a | 2.0 | 22.0 | 257×257 | 748 | 1320 | 0.6 | 49.7 | 3365.8 | Connection D |
| CS-D-2-0.6b | 2.0 | 22.0 | 257×257 | 748 | 1320 | 0.6 | 49.7 | 3412.9 | Connection D |
| CS-D-2-1.0 | 2.0 | 22.0 | 257×257 | 748 | 2118 | 1.0 | 49.7 | 2850.9 | Connection D |

where the $t_{cs}$, $t_g$, $L$, $N_0$, $N_{u,RC}$, $f_c$ and $N_u$, are respectively the thickness of corrugated steel plates, thickness of sandwich grout, length of column, preloading in the RC section, load bearing capacity of RC section, average compressive concrete strength, and ultimate load of specimens; "CS" represents the specimen confined by corrugated steel jacket; the connection type is shown in **Fig. 4**. For example, **CS-A-2-0.4** means that the corner connection of the corrugated steel is the type-A, with thickness of 2 mm and under preload of 0.4$N_{u,RC}$.

The corrugated steel jackets of the specimens are made by splicing several corrugated steel plates. To facilitate the assembly of the corrugated steel plates, four different connection types of the corrugated steel plates were designed, as shown in **Fig. 4**.

(1) **Connection-A**: Angle steel is welded to the corrugated steel plate, and two L-shaped units are bolted together at the corner to form the sleeve.

(2) **Connection-B**: The corrugated steel plates are pre-welded into U-shaped units using angle steel, and then two U-shaped units are connected into a sleeve using a group of bolts.

(3) **Connection-C**: Similar to Connection-A, but the angle steel is replaced with a τ-shaped special connector to increase the flow space for grouting at the corner.

(4) **Connection-D**: Each corrugated steel plate is directly welded to the angle steel in the corner.

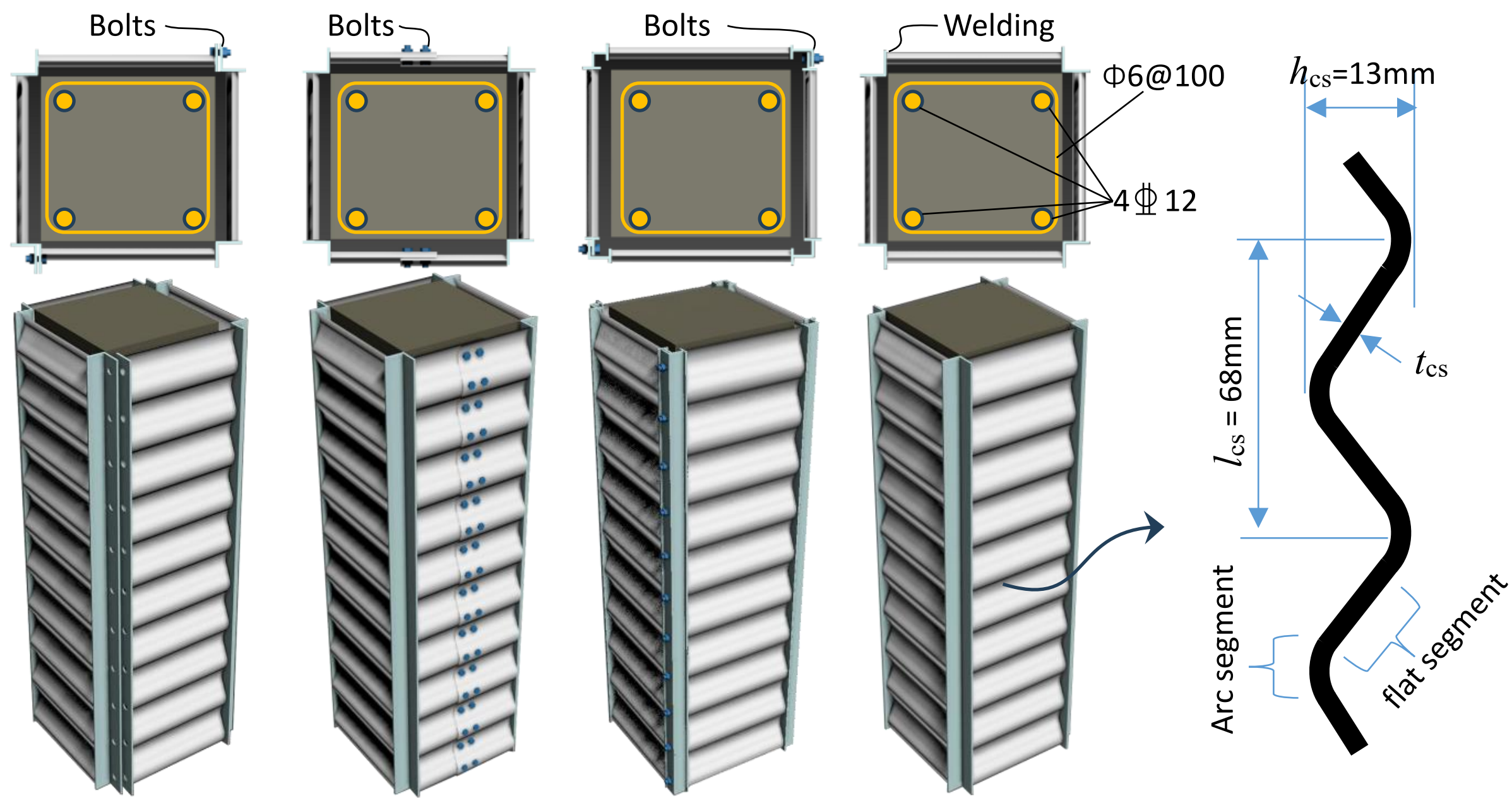

(a) Connection A (b) Connection B (c) Connection C (d) Connection D (e) Corrugation profile

**Fig. 4 Connection types of corrugated steel plates.**

## 2.2 Arrangement of specimens

The workflow of fabrication process of axially compressed short columns and strengthening specimens includes reinforcement cage fabrication, concrete pouring and curing, angle steel welding, corrugated steel plate cutting and welding, strengthening jacket fabrication, applying and removing preload, sleeve enclosure, and interlayer grout pouring.

As shown in **Fig.5**, the jacket is fabricated and installed around the damaged reinforced concrete specimen. Grout is then poured into the interlayer between the concrete and corrugated steel plates. To ensure full filling of the corrugated steel-concrete interlayer, the jacket is fully enclosed, and a two-stage pouring process is adopted:

**(1) First-stage pouring**: The column stands vertically positioned on the ground. The grout is poured from the top to bottom through gaps near the upper end and vibrated to ensure compaction. The initial grout was cured for one day to facilitate the second pouring.

**(2) Second-stage pouring**: The specimen is laid horizontally for secondary pouring. Steel plates are welded to cover the casting window and provide effective confinement to the ends of specimens, avoiding ending local failure.

Upon achieving sufficient strength in the strengthening layer grout, the second-phase loading was initiated until complete failure of specimen occurs.

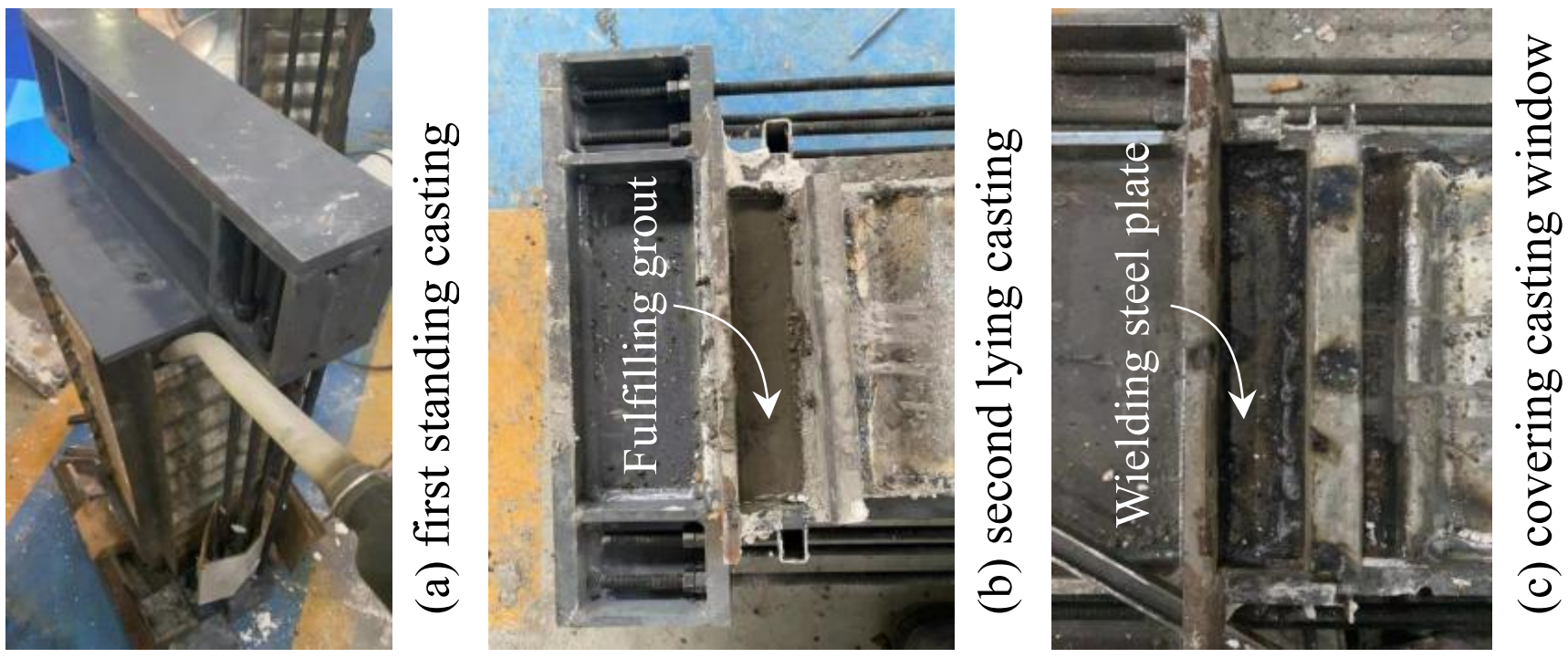

**Fig. 5 Casting grout in enhancing layer.**

### 2.3 Loading and measurement

The testing was conducted in the *Structural Laboratory at Chongqing Jiaotong University*. A four-column servo press machine with a capacity of 10,000 kN was employed to act vertical force to the specimens (see **Fig. 6**). There are two ball hinges locating in the ends of columns. To minimize strain lag effects, low-rate loading was applied, with the entire test conducted at a constant displacement rate of 0.02 mm/min, ensuring that the load is static as closely as possible [42]. Loading was discontinued once the applied load decreased and reached 85% of the ultimate capacity [33]. The preload force of the specimens was measured using force sensors located on the loading plate beneath the press machine. In the reloading process, the axial displacement was monitored using linear variable differential transformers (LVDTs), with data recorded by a static testing system, as shown in **Fig. 6**. Due to the significant variations in strain at the peaks, troughs, and flat segment of the corrugated steel plates, six strain measurement points were installed at these locations on 2 lateral sides of the column at mid-height of the strengthened column specimens. To minimize measurement errors due to bending, small strain gauges (3 mm × 2 mm) were employed on the outer surface of the jackets. Considering that the principal stresses are oriented orthogonally, strain gauges for each measurement point were affixed along both the transverse and longitudinal directions of the corrugation pattern.

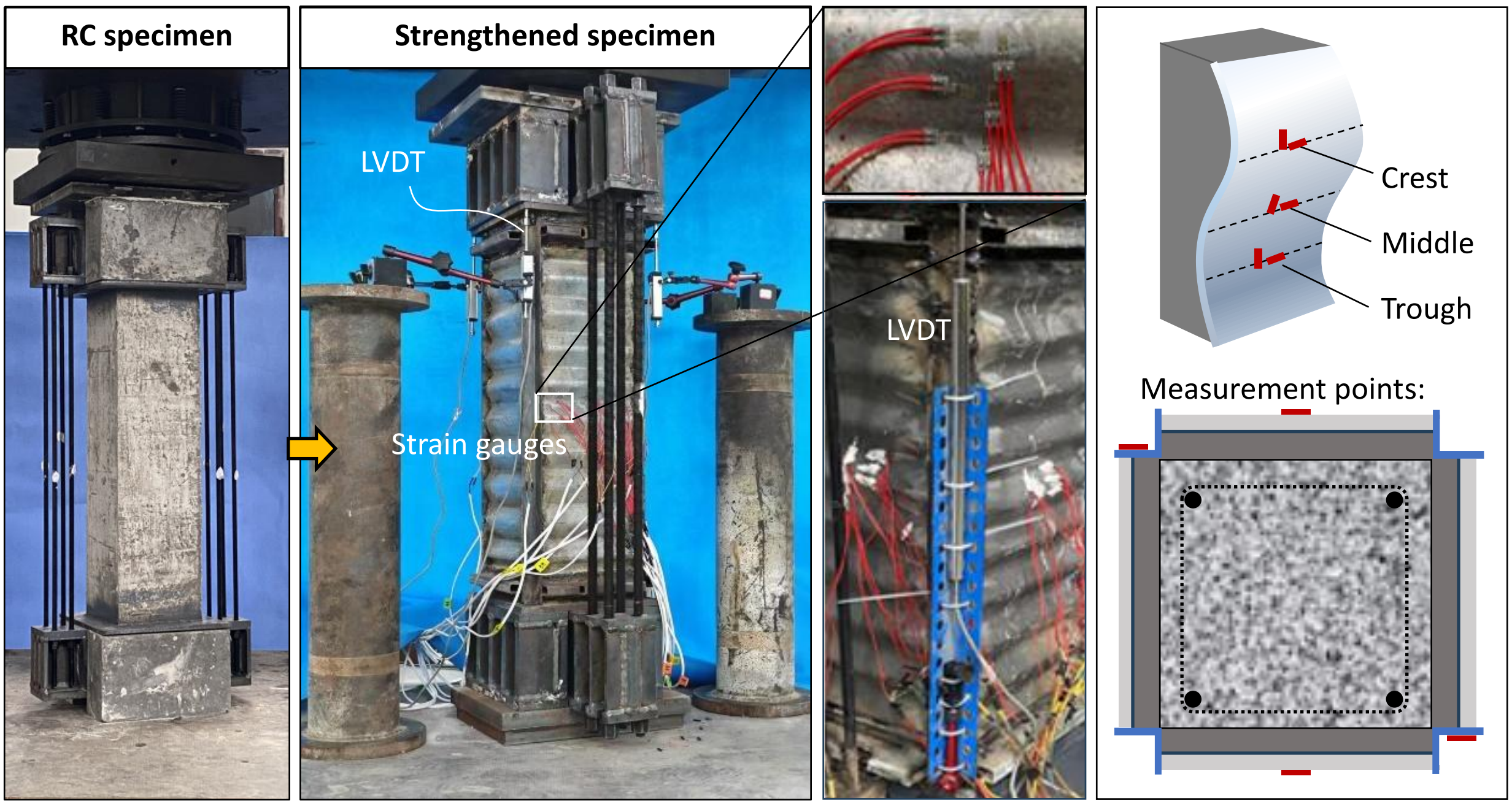

**Fig. 6 Loading and measurement sets.**

## 2.4 Material strength

### 2.4.1 Steel

To test the mechanical properties of steel, samples of corrugated steel plates, angle steels, and reinforcing bars were taken according to the provisions of *Metal Materials Tensile Test – Part 1: Method for Room Temperature Tests* (GB/T 228.1-2010) [43], *Steel for Reinforced Concrete – Part 2: Hot-Finished Ribbed Bars* (GB/T 1499.2-2018) [44], and the standard test methods and definitions for steel mechanical properties (ASTM A370-17a [45]). As shown in **Fig. 7**, due to the varying degrees of cold bending at different local positions on the corrugated steel plates, samples were taken from locations including peaks and flat segment of corrugated steel plates with different thicknesses. During tensile testing, loading was applied at a strain rate of $2.5 \times 10^{-4}$ /s.

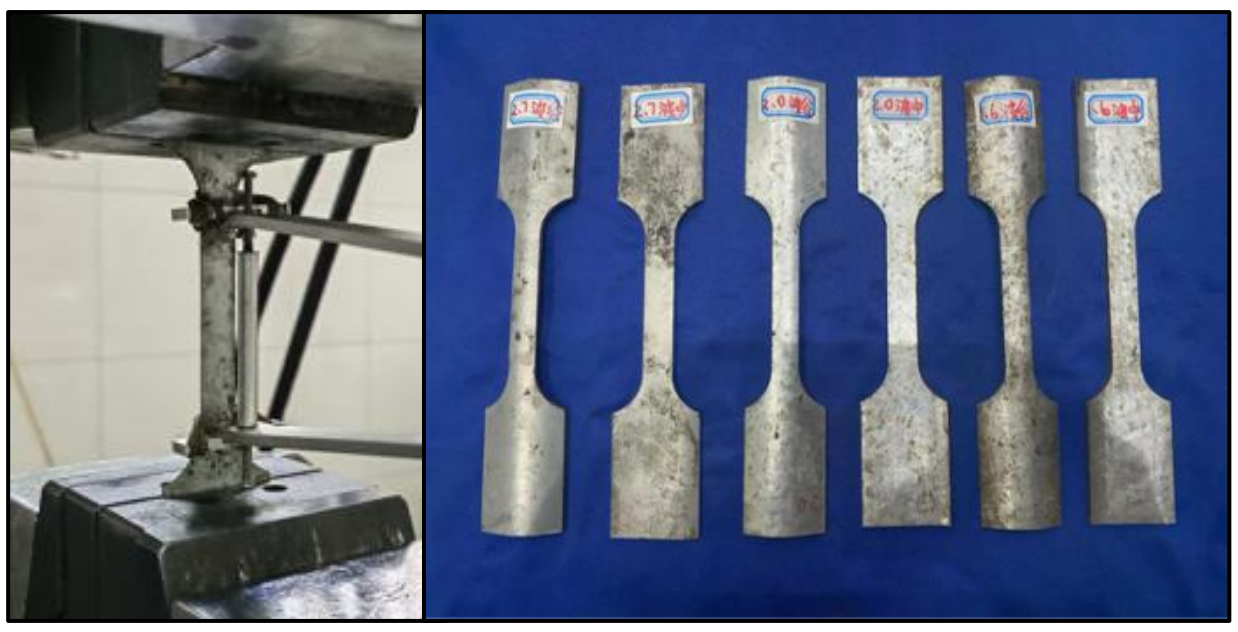

(a) tensile coupons

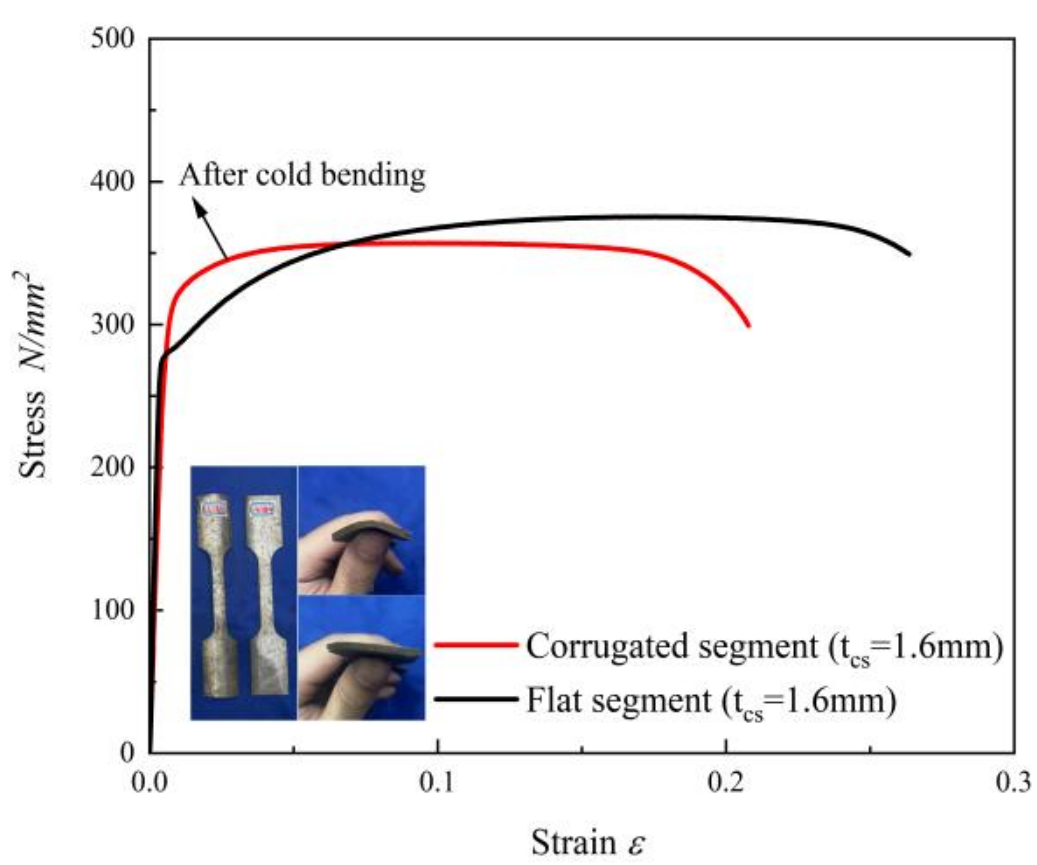

(b) thickness of 1.6mm

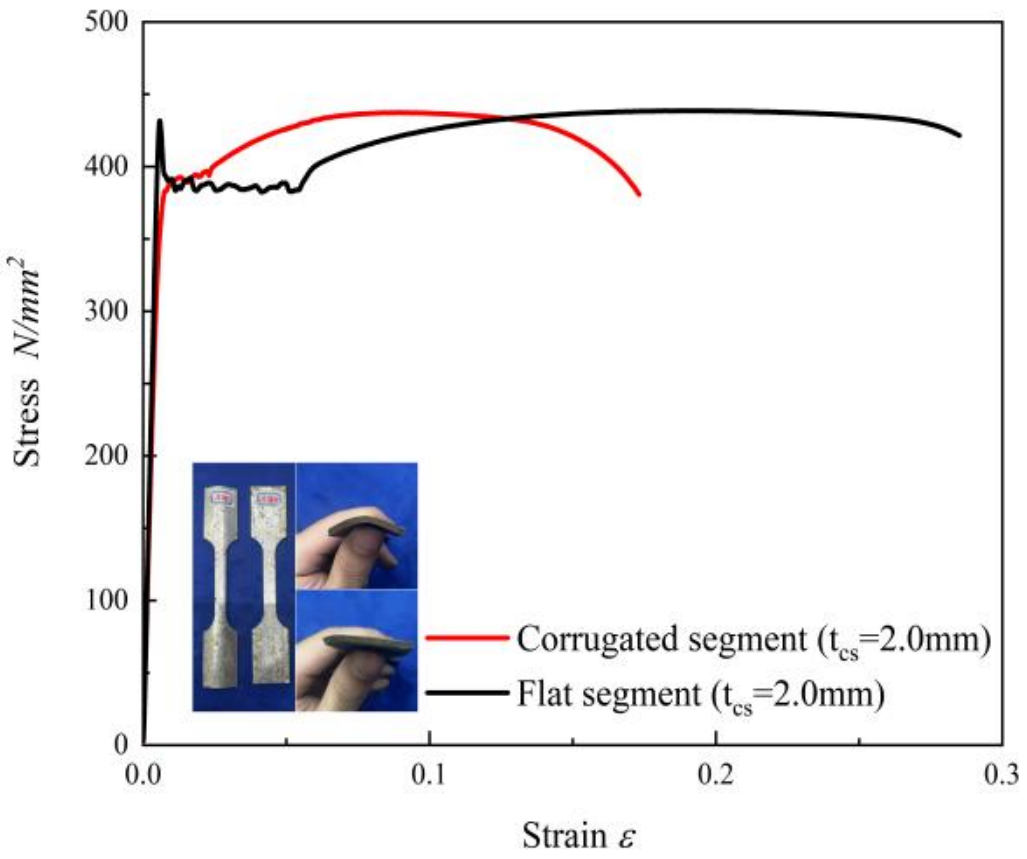

(c) thickness of 2.0mm

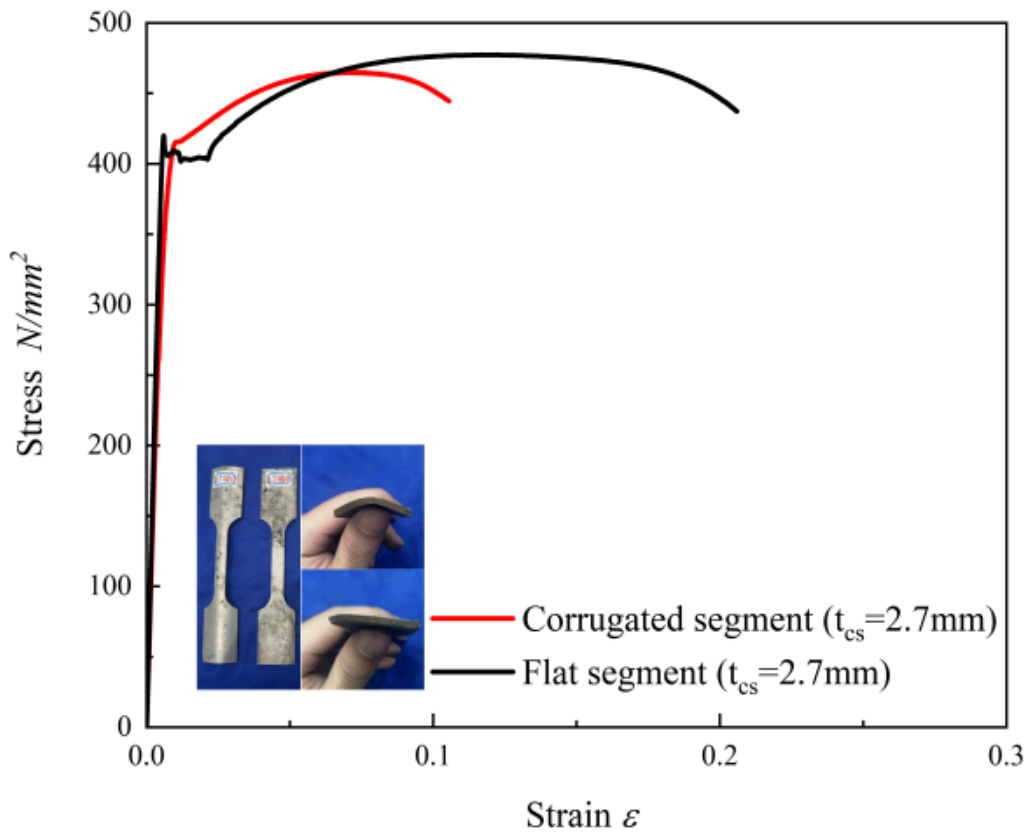

(d) thickness of 2.7mm

**Fig. 7. Behaviour of corrugated steel coupon.**

Note: Due to the varying degrees of cold bending at different local positions of the corrugated steel plate, the yield plateau of the 1.6-mm thick corrugated steel is not clearly visible at the corrugated segment with cold bending.

In practical engineering applications, the yield strength of corrugated steel plates is generally directly taken as the strength of the plate parent material. However, the increase in steel strength due to cold bending is considered as a safety margin. To accurately obtain the actual strength, this study adopted the yield strength after cold bending of the corrugated steel plate. The stress-strain tensile curve of cold-bent corrugated steel plates does not show an obvious yield plateau; therefore, the yield strength was determined using the 0.2% residual strain method. As shown in **Table 2**, the yield strength of corrugated segment is greater than that of flat segment, which matches with the previous studies [33],[35],[36]. Moreover, it can be observed that the corrugated steel plate with the smallest wall thickness (1.6mm) was most affected by cold-forming, resulting in an 18.3% increase in its yield strength. The reduced ductility at the corrugated zones can be observed, which is visible from the earlier drop in the stress-strain curves.

**Table 2. Steel strength.**

| Tensile coupons | $f_y$ /MPa (cold bending contribution) | $f_u$ /MPa (cold bending contribution) | $E$ MPa |
|---|---|---|---|
| Corrugated segment ($t_{cs}$=1.6mm) | 332.81 (+18.3%) | 356.97 (-5.0%) | 198,161 |
| Flat segment ($t_{cs}$=1.6mm) | 281.34 | 375.36 | 201,556 |
| Corrugated segment ($t_{cs}$=2.0mm) | 389.42 (+1.8%) | 437.31 (-0.3%) | 197,212 |
| Flat segment ($t_{cs}$=2.0mm) | 382.35 | 438.60 | 199,214 |
| Corrugated segment ($t_{cs}$=2.7mm) | 413.99 (+3.1%) | 464.85 (-2.6%) | 196,739 |
| Flat segment ($t_{cs}$=2.7mm) | 401.71 | 477.30 | 202,005 |
| Angle steel ($t_{cs}$=4.0mm) | 346.56 | 523.27 | 201,411 |
| Longitudinal bars ($d$=12.0mm) | 443.20 | 589.46 | 199,783 |
| Stirrups ($d$=6.0mm) | 355.20 | 422.68 | 196,358 |

### 2.4.2 Concrete

The test used artificially configured concrete, with the 42.5 ordinary Portland cement; the mass ratio per cubic meter of concrete was as follows: cement/sand/aggregates/water = 376/704/1148/184. The water reducer dosage was 1.2% of the weight of the cementitious materials. During concrete pouring, seven 150-mm concrete cubes were arranged. These test cubes and specimens were all cured under the standard conditions as the specimens for more than 28 days, with an actual loading age of five months. In accordance with the *Code for Test Methods of Physical and Mechanical Properties of Concrete* (GB/T50081-2019) [46], the compressive strength tests were conducted using a 200-ton hydraulic testing machine at the Structural Testing Laboratory of the College of Civil Engineering, Chongqing Jiaotong University. The measured strength data deviating by more than 15% from the average strength value were excluded. In the time period of loading specimens, the average compressive strength of the concrete cubes is 65.39 MPa, and the average compressive strength of prisms [46] is calculated as 49.7 MPa.

### 2.4.3 Grout

According to the *Technical Specifications for Application of Cement-Based Grout Materials* (GB50488-2015) [47], when the maximum aggregate size in cement-based grout materials does not exceed 4.75 mm, the prism with dimensions of 40 mm×40 mm×160 mm is recommended to test the compressive strength. Considering the thin reinforcement layer thickness ($t_g$=20 mm), six prisms with dimensions of 40 mm × 40 mm × 160 mm was used. The compressive strength test was conducted in accordance with the *Cement Mortar Strength Testing Method (ISO Method)* (GB/T 17671-2021) [48]. The average compressive strength of the grout prism was 37.1 MPa.

## 3. Experiment Results

### 3.1 Failure modes

All specimens showed no significant damage until the applied load reached its peak. Only when the load began to decrease, minor buckling occurred at the corner steel angles, while the specimens remained intact overall, with no obvious lateral bulging (**Figs. 8 and 9**). This is because the corrugated profile has a greater moment of inertia than flat plates, allowing it to resist outward deformation more effectively through its transverse stiffness and interlocking mechanism. It indicates that even with a steel plate thickness as thin as 1.6 – 2.7 mm, i.e., width-to-thickness ratios ranging from 74 to 125, sufficient confinement can still be achieved. It can be validated by comparisons of the load-displacement curves obtained from experimental tests and stress analysis results.

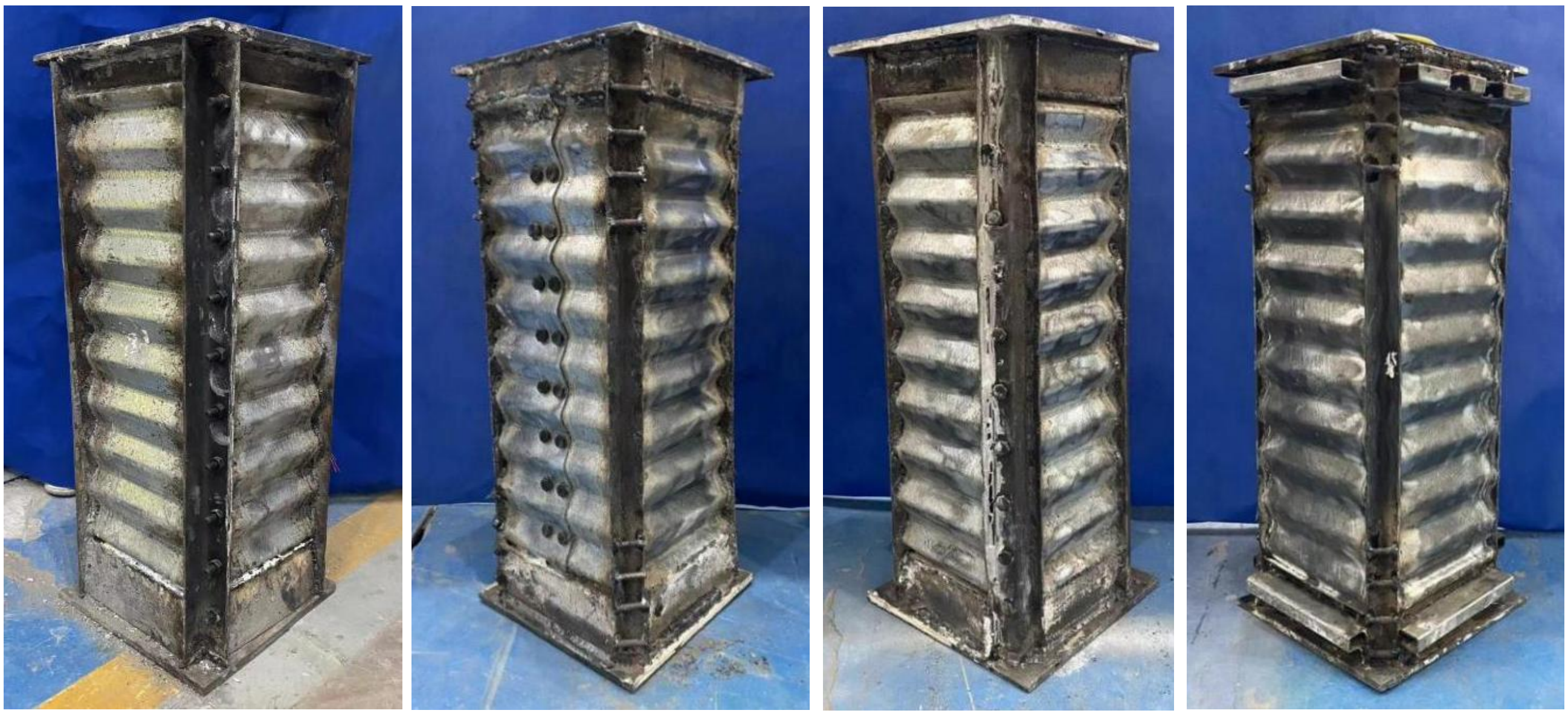

(a) connection A (b) connection B (c) connection C (d) connection D

**Fig. 8. Whole failure of the specimens.**

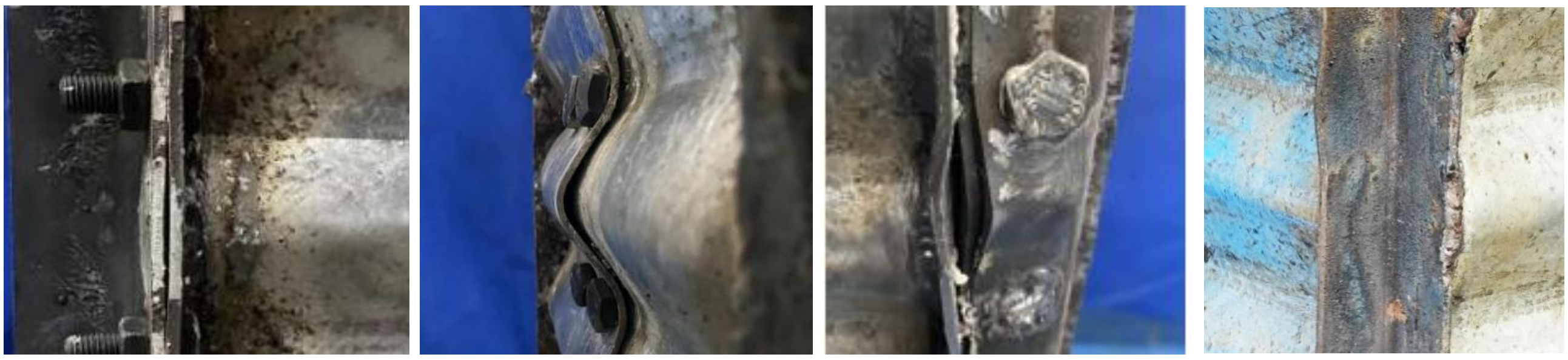

(a) connection A (b) connection B (c) connection C (d) connection D

**Fig. 9. Local failure of the specimens.**

For the specimen with Connection A, no significant damage occurred throughout the entire loading process. At 75% of the ultimate load, slight local buckling and loosening were observed between the two bolts at the mid-section's steel angles and flange plates.

For the specimen with the Connection B, during the initial loading to peak load, no significant changes occurred on the external surface of the steel plates, angles, and the central connection of the corrugated steel sheets. As displacement increased further and loading was applied until the load reached 80% of the ultimate load (i.e., 2500 kN), minor buckling occurred at the corner angles. When the load was reduced to 60% of the ultimate load (1900 kN), a peeling trend was observed between the two corrugated steel sheets at the central connection, with slight lateral bulging of the adjacent steel plates.

For the specimen with the Connection C, no significant changes occurred until the load was reduced to 90% of the ultimate load. At that moment, a muffled sound was heard, suggesting potential crushing of the core concrete or grout between the layers. Continued loading to 50% of the ultimate load resulted in local buckling and separation of the two corner connected parts, while the bolts remained intact without loosening.

As the simpler welding configuration induces less stress concentration than bolted connections, no obvious failure phenomena were observed during loading in specimens with Connection D, including those under different preload levels. It was not until the load decreased to 70% of the ultimate capacity that periodic post-peak buckling occurred along the corrugation. This multi-wave buckling behavior indicates high stability in the thin-walled plate.

### 3.2 Load-Displacement Curves

During the early stages of loading, even though the core concrete had developed some plastic behavior, the rest of the concrete remained in the elastic working stage due to the enhanced section. As a result, the curve initially exhibited roughly linear-elastic behavior. When the load increased and exceeded the original ultimate bearing capacity of the RC section, the specimen curve began to show curvature, with the onset of plastic deformation. Upon reaching the limit state, the load started to decrease, and the peak deformation corresponding to the load was approximately twice the original RC section’s peak deformation. The ductility in the descending segment was also significantly better than that of the plain RC specimens.

From the comparison, it can be observed that after being wrapped with sleeves and having its section enlarged, the ultimate bearing capacity was significantly strengthened. Although some cracking occurred under the initial axial compression, potentially reducing the strength, the constrained concrete’s strength was enhanced due to the confinement effect. As a result, the ultimate bearing capacity increased from 2117.7 kN to 2850.9 - 3543.4 kN, representing an increase of 34.6% to 67.3%.

The influence of thickness of corrugated steel jacket is shown in **Fig. 10**. There is no obvious impact on the load-displacement (*N-Δ*) curves of specimens with different thickness of the corrugated steel jacket. It suggests that the confinement of the corrugated steel jacket may not be taken full advantage. As well, the enhanced ultimate capacity and ductility are both far greater than that of the original RC section, indicating that it may be unnecessary to design a thick jacket.

The effect of different connection types can be compared in **Fig. 11**. It shows little difference among the four connection types in terms of ultimate bearing capacity. However, regarding connection strength, the specimen Connection D with full welding at the corner results in the most reliable radial confinement. As a result, its ultimate bearing capacity (3361.3 kN) is 9.8% greater than the mean value other 3 specimens. In specimen Connection C, a sharp drop in load occurred when residual load was reduced to 3000 kN, likely due to the steel plate at the corner experiencing reduced effective cross-section and decreased compressive stiffness because of openings or dense welds, leading to local buckling. However, after reaching around 2000 kN, the load decline slowed down as residual capacity was regained due to circumferential confinement. For specimen Connection B, its bearing capacity is moderate. After loading ended, only slight gaps were observed at the wave plate connections. The specimen Connection A exhibits unique behaviour: its load bearing capacity remained stable for a period after reaching peak load. However, when compressive displacement reached 4 mm (compression strain = 0.005), local instability failure occurred due to multi-wave buckling of the angle steel, causing a sudden drop in bearing capacity. As lateral confinement became effective, the rate of load decline slowed down.

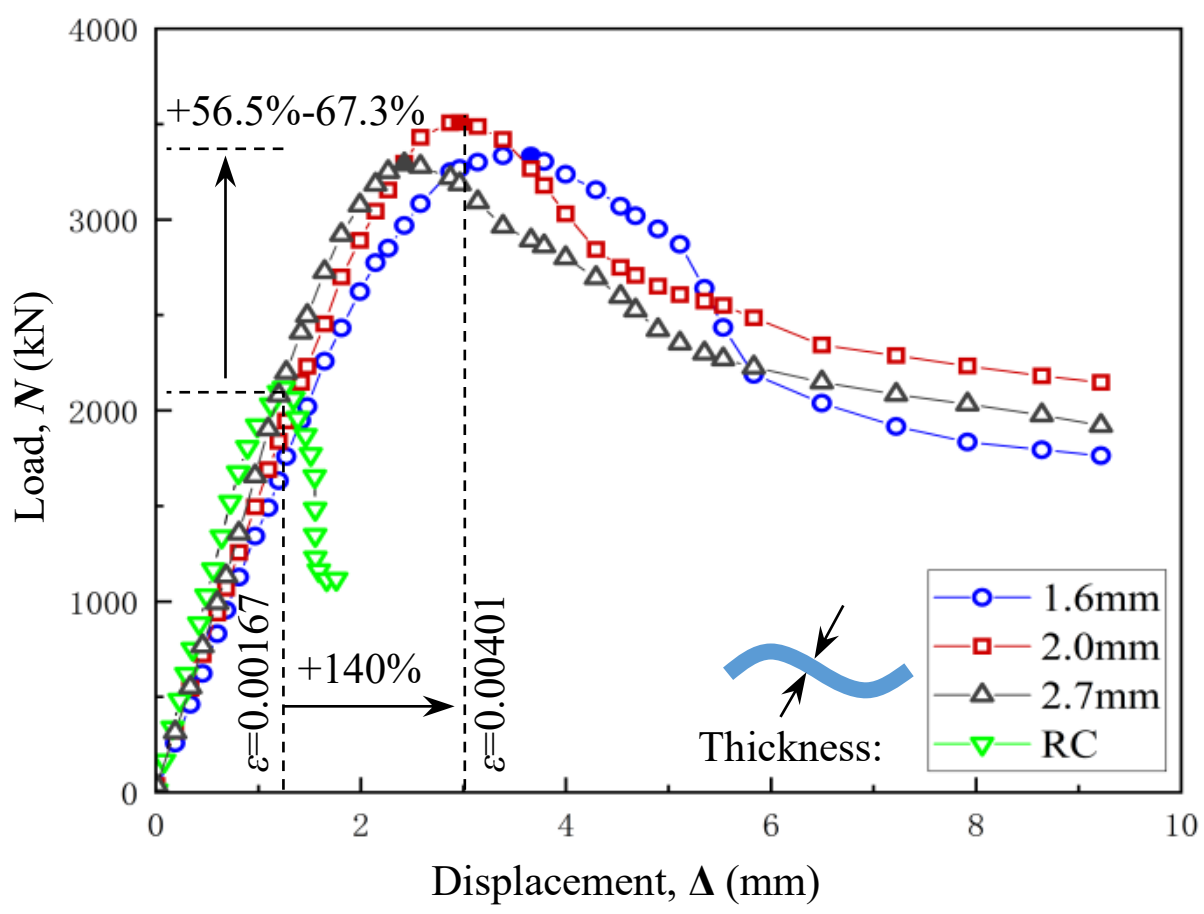

**Fig. 10. Influence of thickness of jacket.**

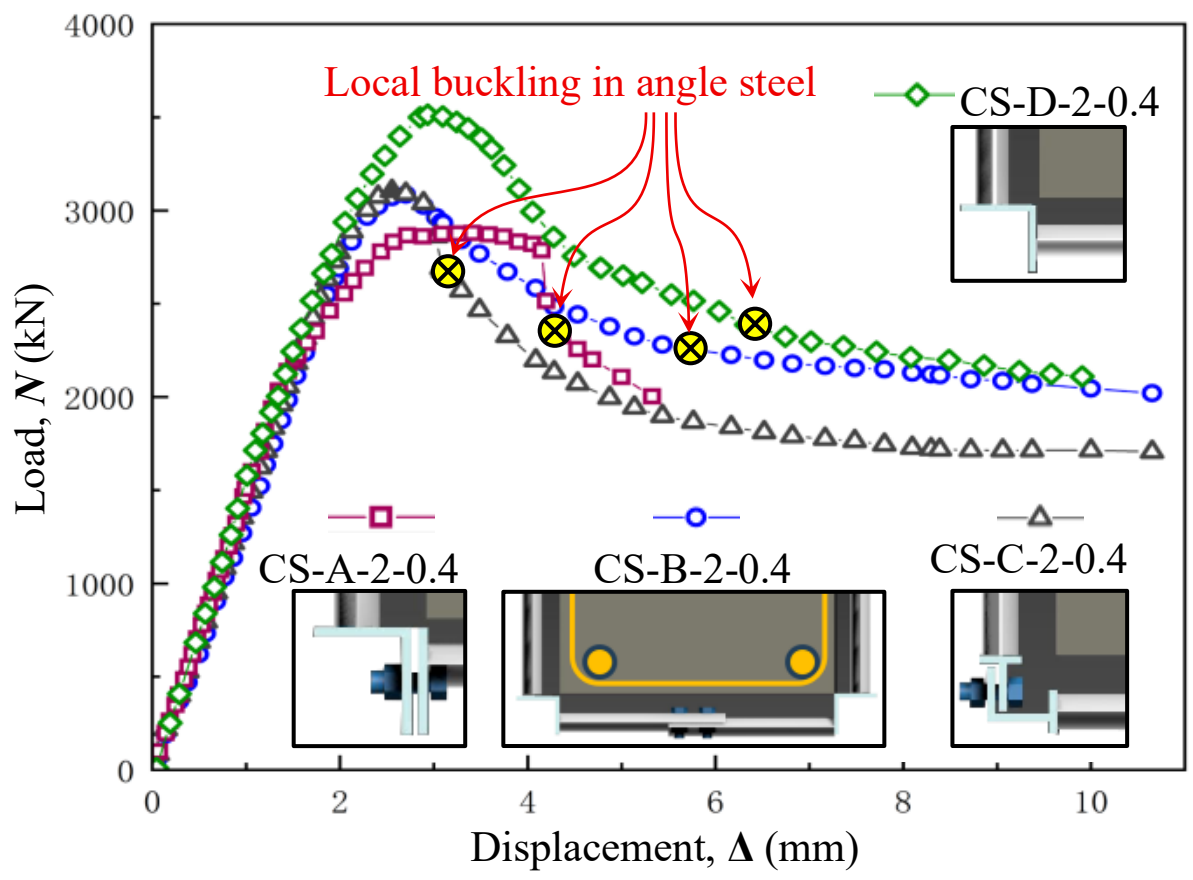

**Fig. 11. Influence of connection type.**

As shown in **Fig. 12**, the influence of preloading level on the performance enhancement of strengthened RC sections is also evident in the load–displacement ($N$–$\Delta$) curves. The ultimate load capacity of the strengthened specimens decreases with increasing preload level. Specifically, a higher initial preload reduces the strengthening effect. For specimens preloaded between 0.4 and 0.6 of $N_{u,RC}$, the increase in bearing capacity relative to the unreinforced RC section declines from 67.3% to 60.0%. In the extreme case where the preload reaches $N_{u,RC}$, the concrete softens significantly, accompanied by a sharp reduction in residual strength. Nevertheless, the ultimate bearing capacity still shows a 34.6% increase compared to the original RC section. This softening influences not only load capacity but also deformation behavior. Despite initial damage in the concrete, the jacketed RC column effectively behaves as a new steel–concrete composite member, analogous to a concrete-filled steel tubular column whose ductility is known to correlate with the confinement factor. Due to the reduced concrete strength in the inner RC section resulting from preloading, a higher confinement factor is achieved, as established in previous studies [33],[49],[51][52][53][54]. Consequently, the specimen preloaded to 100% $N_{u,RC}$ exhibits greater peak deformation than those with lower preload levels (0.4 and 0.6 of $N_{u,RC}$).

Briefly, in a jacketed column after unloading, greater pre-damage causes concrete softening and increased lateral expansion under re-compression. This dilation promotes a tighter interaction between the core concrete and the outer jacket, activating stronger passive confinement and ultimately enhancing ductility. The findings are also supported by subsequent ductility index analysis (**Fig.13(c)** and **Table. 3**) and high level of horizontal stress of CS-D-2-1.0 (**Fig.14 (a)**).

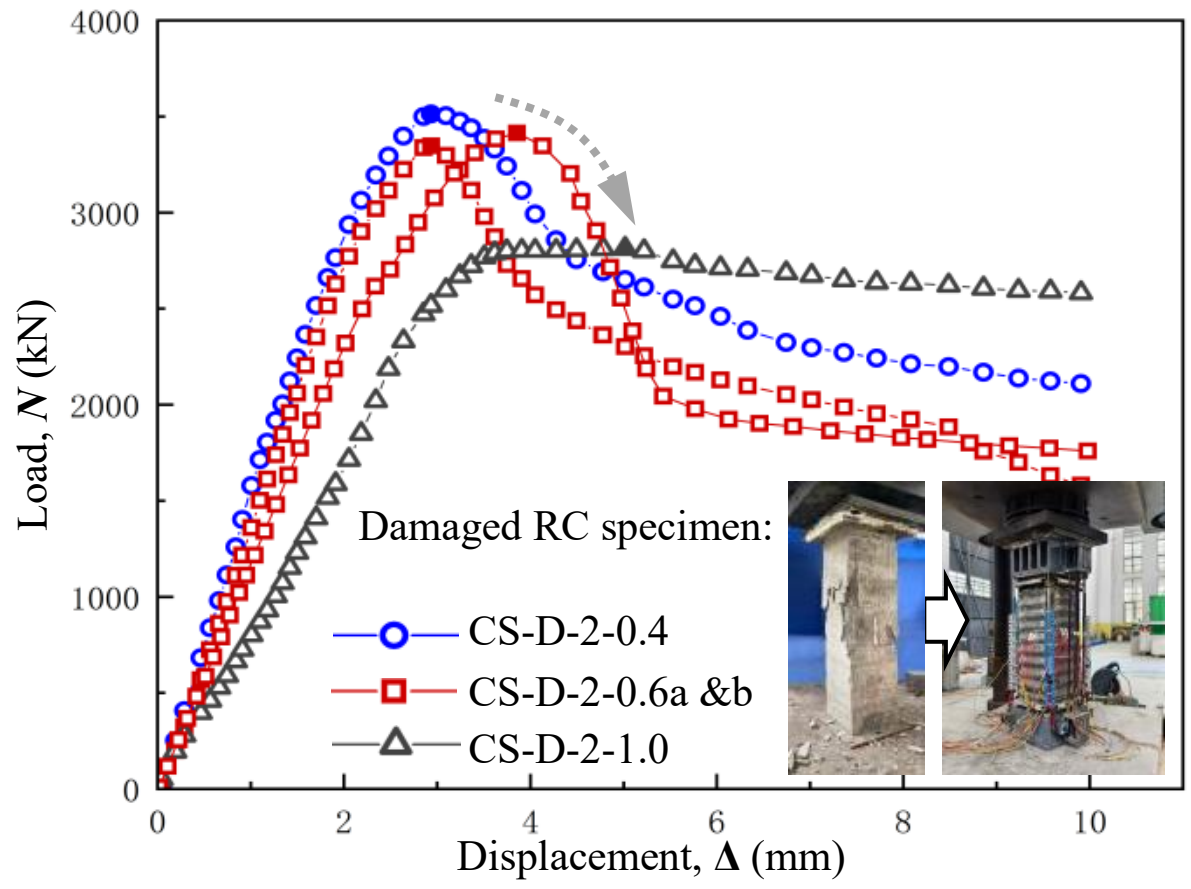

**Fig. 12. Influence of initial damage.**

### 3.3 Strength index SI and ductility index DI

As the concrete strength was not completely the same in the test, a further analysis is necessary based on dimensionless indicators. The strength index (*SI*) [49] can be used for estimate the load carry capacity of an ordinary reinforced concrete section strengthened by corrugated steel jacket, which can be calculated as Eq. (1). Because the corrugated steel plates with a sine-curve corrugation may not directly carry the axial load[75], it is not taken into this calculation.

$$SI = \frac{N_{\mathrm{u}}}{f_{\mathrm{c}}A_{\mathrm{c}} + f_{\mathrm{c,g}}A_{\mathrm{g}} + f_{\mathrm{y,s}}A_{\mathrm{s}} + f_{\mathrm{y,L}}A_{\mathrm{L}}} \tag{1}$$

where $N_{\mathrm{u}}$, $f_{\mathrm{c}}$, $f_{\mathrm{c,g}}$, $f_{\mathrm{y,s}}$, $f_{\mathrm{y,L}}$, $A_{\mathrm{c}}$, $A_{\mathrm{g}}$, $A_{\mathrm{s}}$ and $A_{\mathrm{L}}$ are respectively the ultimate bearing capacity of the column, the compressive strength of the concrete and grouting material, the yield strength of the reinforcing bars and angle steel, and their corresponding cross-sectional areas.

The ductility index, *DI*, well represents the deformation capacity and energy absorption of a structural member [49]. As shown below:

$$DI = \varepsilon_{85\%} / \varepsilon_{\mathrm{y}} \tag{2}$$

where $\varepsilon_{85\%}$ and $\varepsilon_{\mathrm{y}}$ are the axial compressive strain corresponding to $0.85N_{\mathrm{u}}$ in the descending stage and yield compressive strain. The $\varepsilon_{\mathrm{y}}$ can be obtained using the farthest point method proposed by Feng P et. al. (2017) [50].

The confinement factor, $\xi$, is an important index to assess the mechanical behaviour of concrete-filled tubular columns [49]. If the columns were not damaged, the $\xi$ should be slightly minor. Because the confining tube is constructed after the damaged of concrete, the residual compressive strength of the concrete with damage ($f_c'$) should be employed to calculate $\xi$:

$$\xi = \frac{f_{\mathrm{y,cs}}A_{\mathrm{cs}} + f_{\mathrm{y,L}}A_{\mathrm{L}}}{f_{\mathrm{c}}'A_{\mathrm{c}} + f_{\mathrm{c,g}}A_{\mathrm{g}}} \tag{3}$$

As shown in **Fig. 13**, the strength index (*SI*) and ductility index (*DI*) exhibit a strong dependence on two key parameters: the confinement factor ($\xi$) and the preload ratio ($N_0/N_{\mathrm{u,RC}}$). Increased $\xi$ enhances both *SI* and *DI*, reflecting enhanced confinement effects. In contrast, the $N_0/N_{\mathrm{u,RC}}$ influences *SI* and *DI* differently. The increase in $N_0/N_{\mathrm{u,RC}}$ reduces *SI* due to more cracks developing in the concrete, while it simultaneously enhances *DI* due to the softening behavior of the concrete. In terms of the connection type of corrugated steel jacketing, CS-B-2-0.4 shows a similar *DI* to CS-D-2-0.4, indicating that the bolts connection may have a satisfactory strength as the welding in the corner if the design of connection strength is rational.

**Table 3. *SI* and *DI* of specimens.**

| No. | Specimens | Preload ratio $N_0/N_{u,RC}$ | $\alpha$ | $\xi$ | $N_{u0}$ (kN) | $N_u$ (kN) | *SI* | *DI* |
|---|---|---|---|---|---|---|---|---|
| 1 | RC | 0 | 0 | 0 | / | 2117.7 | 1.00 | 1.86 |
| 2 | CS-A-2-0.4 | 0.4 | 5.45% | 0.43 | 2956.3 | 2957.4 | 1.00 | 2.95 |
| 3 | CS-B-2-0.4 | 0.4 | 4.25% | 0.34 | 3068.0 | 3101.1 | 1.01 | 2.21 |
| 4 | CS-C-2-0.4 | 0.4 | 4.88% | 0.34 | 2989.9 | 3125.6 | 1.05 | 1.70 |
| 5 | CS-D-1.6-0.4 | 0.4 | 3.64% | 0.25 | 3082.9 | 3361.3 | 1.09 | 2.59 |
| 6 | CS-D-2.7-0.4 | 0.4 | 5.23% | 0.44 | 3082.9 | 3314.0 | 1.07 | 2.43 |
| 7 | CS-D-2-0.4 | 0.4 | 4.22% | 0.34 | 3082.9 | 3543.4 | 1.15 | 2.25 |
| 8 | CS-D-2-0.6a | 0.6 | 4.22% | 0.34 | 3052.1 | 3365.8 | 1.10 | 1.91 |
| 9 | CS-D-2-0.6b | 0.6 | 4.22% | 0.34 | 3052.1 | 3412.9 | 1.12 | 2.14 |
| 10 | CS-D-2-1.0 | 1.0 | 4.22% | 0.35 | 2942.0 | 2850.9 | 0.97 | 3.51 |

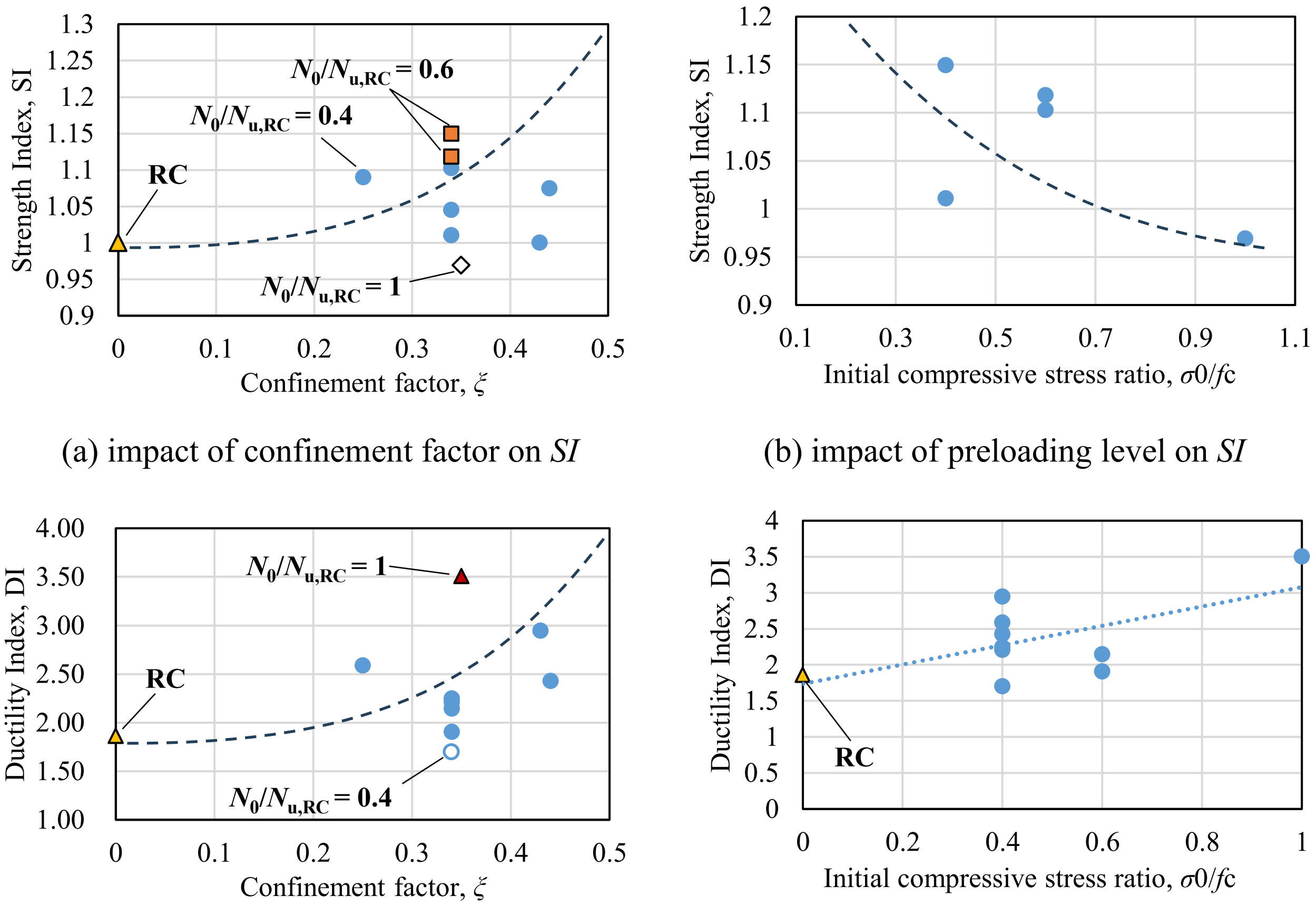

(a) impact of confinement factor on *SI* (b) impact of preloading level on *SI*

(c) impact of confinement factor on *DI* (d) impact of preloading level on *DI*

**Fig. 13. *SI* and *DI*.**

## 4. Confinement Mechanism

### 4.1 Load Carrying Mode

Previous studies on corrugated steel concrete columns have shown that the primary role of the corrugated steel is to provide circumferential confinement [32,33,35]. Similarly, in reinforced concrete sections strengthened with corrugated steel plates, the outer corrugated steel plates bear minimal longitudinal loads and mainly provide lateral confinement. The parts directly carry the axial load include the original RC section, the reinforcement layer grout, and the corner connection components of the sleeves. It can be verified by examining the development patterns of strain and stress.

#### 4.1.1 Strain of corrugated jacket

The strain of corrugated steel plates is influenced by their corrugation profile, with variations in strain occurring across the crests, troughs, and flat segment within one period of corrugation. For transverse strain, overall tensile strain is relatively small, possibly due to the placement of strain gauges farther from the corners, where the confinement is the greatest in a square or rectangular cross section [51-54]. In the middle regions, steel does not need to provide significant confinement, resulting in reduced transverse tensile strain. For longitudinal strain, the trough of the corrugated steel plate exhibits a more pronounced compressive longitudinal strain, indicating that the concave trough may take on part of the load. Conversely, the crest shows nearly zero longitudinal strain, suggesting that the convex crest relieves some of the applied stress. This phenomenon, known as the "bending effect," was proposed by Wang and Yang et al. (2019) [33] under compressive loading conditions.

#### 4.1.2 Stresses of corrugated jacket

To investigate the load-bearing behavior of the corrugated jacket, elastic-plastic analysis [32,33,55] was employed to calculate the stresses in the corrugated steel. The stress components (horizontal and vertical stresses) as well as the Mises stress were determined at the crest, flat segment (middle), and trough on the external surface of the corrugated plates, based on the measured horizontal and vertical strains.

**(1) Plane elastic-plastic analysis method**

For the stress at different positions on the corrugation, the respective Poisson's ratio, elastic modulus, and steel yield strength are incorporated into the calculations to obtain a more precise analysis of the corrugated steel plate stress. Since the corrugations in the steel plates run horizontally, the principal stresses in the columns under axial compression are primarily in the horizontal and vertical directions. Given that the shear stress within the corrugated steel plates is relatively small, only these two directional normal stresses were analyzed in the test.

In the elastic stage, the stress components meet the generalized Hooke's Law:

$$\begin{bmatrix}\sigma_h\\ \sigma_v\end{bmatrix}=\frac{E_s}{1-\mu_s^2}\begin{bmatrix}1 & \mu_s\\ \mu_s & 1\end{bmatrix}\begin{bmatrix}\varepsilon_h\\ \varepsilon_v\end{bmatrix} \tag{4}$$

where $\varepsilon_h$, $\varepsilon_v$, $\sigma_h$, $\sigma_v$, $E_s$ and $\mu_s$ are respectively the horizontal strain, vertical strain, horizontal stress, vertical stress, elastic modular and Poisson ratio.

In the elastic-plastic stage, the $E_s$ and $\mu_s$ can be replaced by the tangential modular ($E_s^t$) and elastic-plastic Poisson ratio ($\mu_{sp}$)[56]:

$$\begin{bmatrix}d\sigma_h\\ d\sigma_v\end{bmatrix}=\frac{E_s^t}{1-\mu_{sp}^2}\begin{bmatrix}1 & \mu_{sp}\\ \mu_{sp} & 1\end{bmatrix}\begin{bmatrix}d\varepsilon_h\\ d\varepsilon_v\end{bmatrix} \tag{5}$$

$$E_s^t=\frac{(f_y-\sigma)}{(f_y-f_p)}\frac{\sigma}{f_p}E_s \tag{6}$$

$$\mu_{sp}=0.283+0.167\frac{\sigma-f_p}{f_y-f_p} \tag{7}$$

The stress-strain relation is more complex in the plastic stage, which can be expressed by:

$$\begin{bmatrix}d\sigma_h\\ d\sigma_v\end{bmatrix}=\frac{E_s}{Q}\begin{bmatrix}\sigma_v'^2+2p & 2\mu_s p-\sigma_h'\sigma_v'\\ 2\mu_s p-\sigma_h'\sigma_v' & \sigma_v'^2+2p\end{bmatrix}\begin{bmatrix}d\varepsilon_h\\ d\varepsilon_v\end{bmatrix} \tag{8}$$

**(2) Strain and stress analysis result**

Since strain was recorded on both sides of the cross-section at the same position in the corrugation, two load-strain curves were obtained, as shown in **Figs. A-1** and **A-2**.

The horizontal strain or stress of the corrugated steel jacket represents the degree of confinement. During the ascending stage of loading, it exhibits a slight tensile increase. In the descending stage, the horizontal strain generally decreases rapidly at both the crest and trough of the corrugation, or even become compressive. Only the middle of corrugation remains tensile. This behavior may be attributed to the Poisson effect, where increasing vertical compressive strain suppresses horizontal tensile strain. For the specimen CS-D-2-1.0 with complete damage, the horizontal strain at one of the measured points has approached the yield strain, and horizontal stress can reach near 200 MPa at the ultimate load. It indicates that the initially damaged concrete has a more significant lateral expansion, which is prevented by the corrugated steel sheet. However, the horizontal strain and stress become minor with less damage, suggesting that the cracking and softening of the damaged concrete can increase the confinement effect in the round strengthened RC column. The corrugated steel sheet should have exhibited good confinement [32], whereas, it was not taken full advantage due to the reduced confinement development in a square section[84].

For vertical strain and stress, similar behavior can be observed in concrete-filled corrugated steel tubes with circular cross-sections [32]. During compression of the corrugated steel together with the column, the outer-surface fibers are stretched at the peak and flat segments of the corrugation, where tensile and compressive stresses cancel each other. Conversely, the trough is compressed more with the superimposing of column compressive stress, resulting in greater compressive stress. This superposed stress causes the Mises stress at the trough to increase significantly faster than at the crest and flat regions of the corrugation (**Fig. A-2**). As shown in **Fig. 14 (b)**, the average longitudinal stress of all specimens in the corrugated steel was only -3.591 MPa by using weighted mean values[32], indicating that the axial load directly borne by the corrugated steel plate can be negligible. This finding is consistent with previous relevant research conclusions[33,35,36,54].

In general, the corrugated steel plates at different positions exhibit low stress levels, with most horizontal stresses and vertical stresses not exceeding 100 MPa, and the Mises stress not reaching the yield strength value (**Fig. 14 (c)**). The exception occurs when severe plastic strain develops in the RC section with complete damage. Given this low stress level, the thickness of the corrugated steel does not need to be excessively large; it only needs to satisfy the construction requirements for local stability and connection reliability.

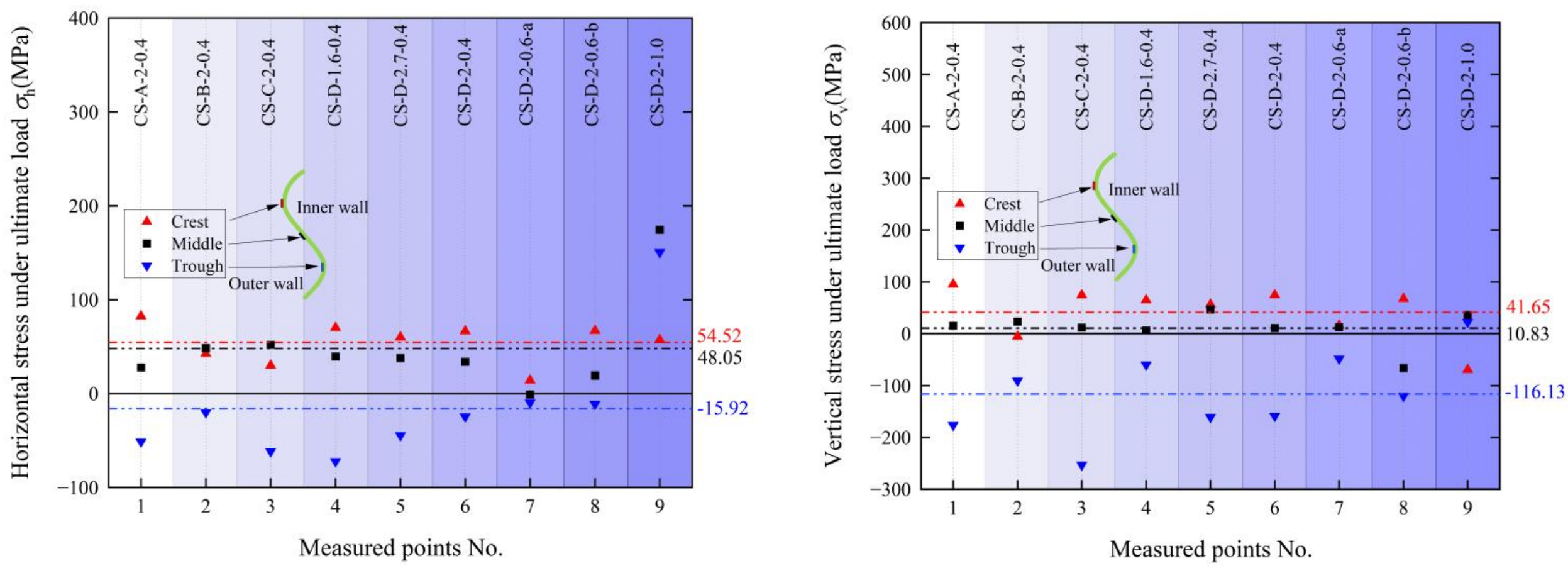

(a) horizontal stresses

(b) vertical stresses

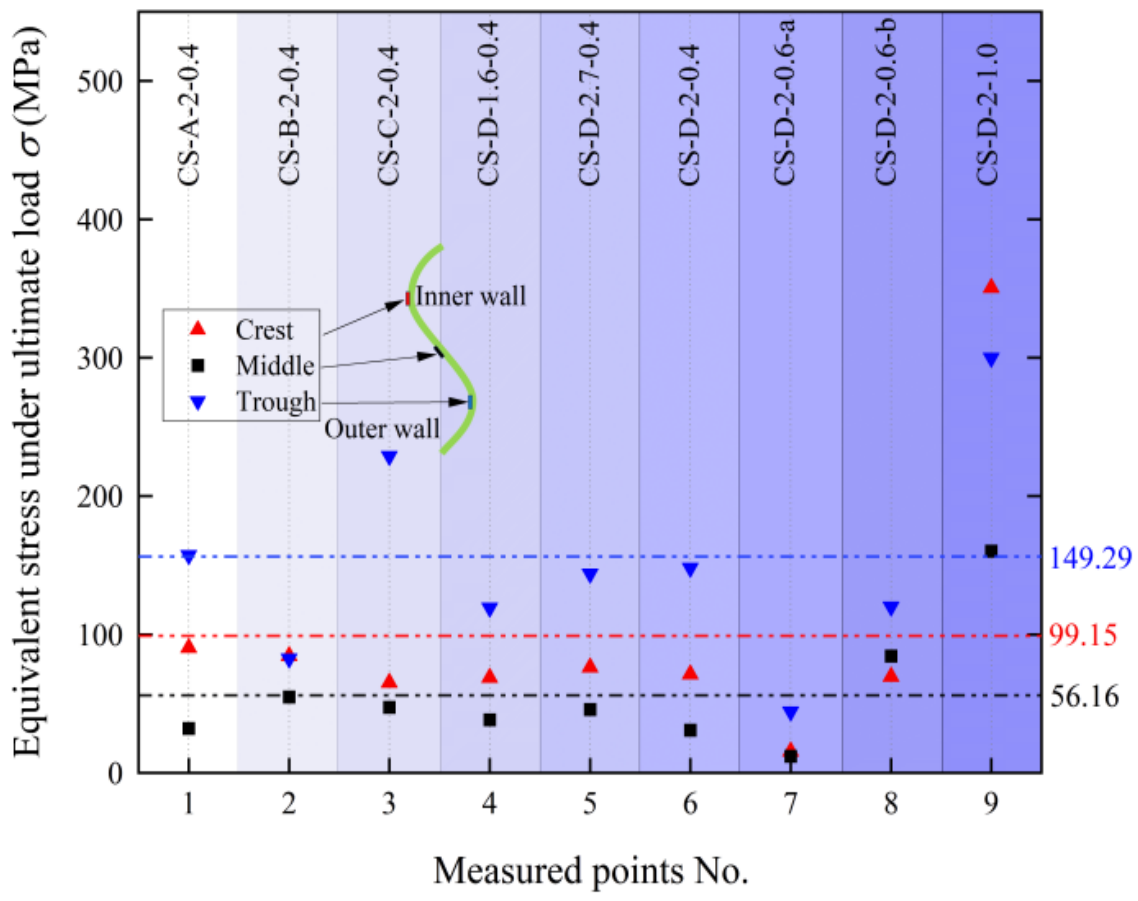

(c) Mises stresses

**Fig. 14. Stress components and Mises stress at ultimate load.**

As shown in **Fig. 15.**, the centrally located strongly confined concrete plays a dominant role in the load carrying (60%-80%), while those of steel components such as reinforcement and angle steel are relatively minor (< 20% in general). For such a strengthening method with a thin sandwich grout layer, the structure primarily relies on the external confinement of the core concrete to generate significant confinement effect, thereby enhancing its load-bearing capacity.

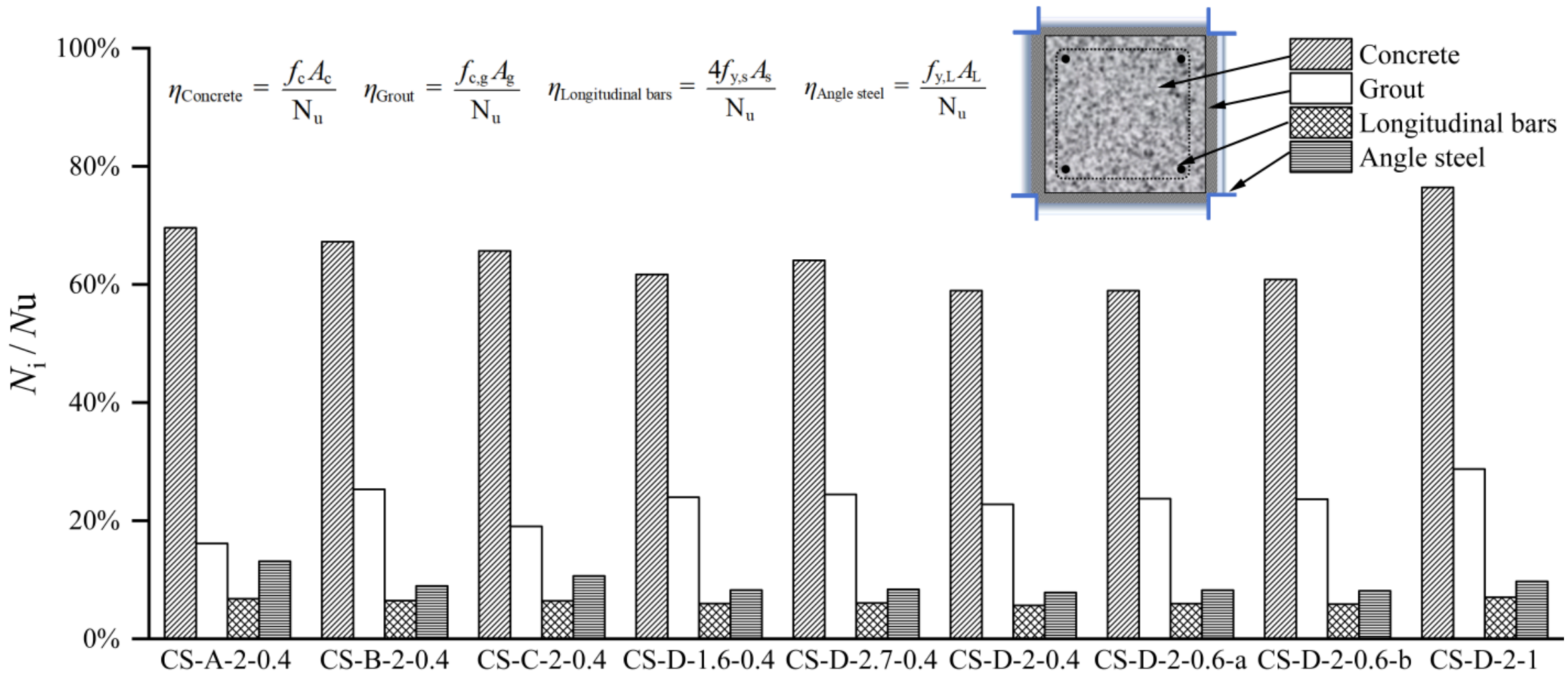

**Fig. 15. The load-bearing proportion of each part**

## 5. Ultimate axial load bearing capacity

### 5.1 Ultimate state mechanics

To facilitate the calculation and design of the ultimate bearing capacity of reinforced concrete sections after strengthening, the stress distribution across the section at the ultimate load-bearing state is analyzed as shown in **Fig. 16**. The reinforced concrete sections after strengthening are no longer ordinary but have become a composite section confined by the jacket, with its load-bearing mechanism similar to that of square steel tube concrete columns [54]. In the regions close to the corners, the confinement is strong, whereas, it becomes weaker near the centerlines on the sides. Therefore, when the column is under compression, the concrete and grout in regions A1 and A3 (confined areas) are subjected to triaxial compression. Conversely, in the areas A2 and A4 that is weakly confined or unconfined, the concrete may even experience uniaxial compression. To accurately calculate the ultimate bearing capacity of such sections, it is necessary to analyze each region separately and then combine the results. In the calculations, the boundary between the strongly confined and weakly confined zones is assumed to be a quarter-circle arc. Based on this assumption, the areas of each partition can be computed as detailed.

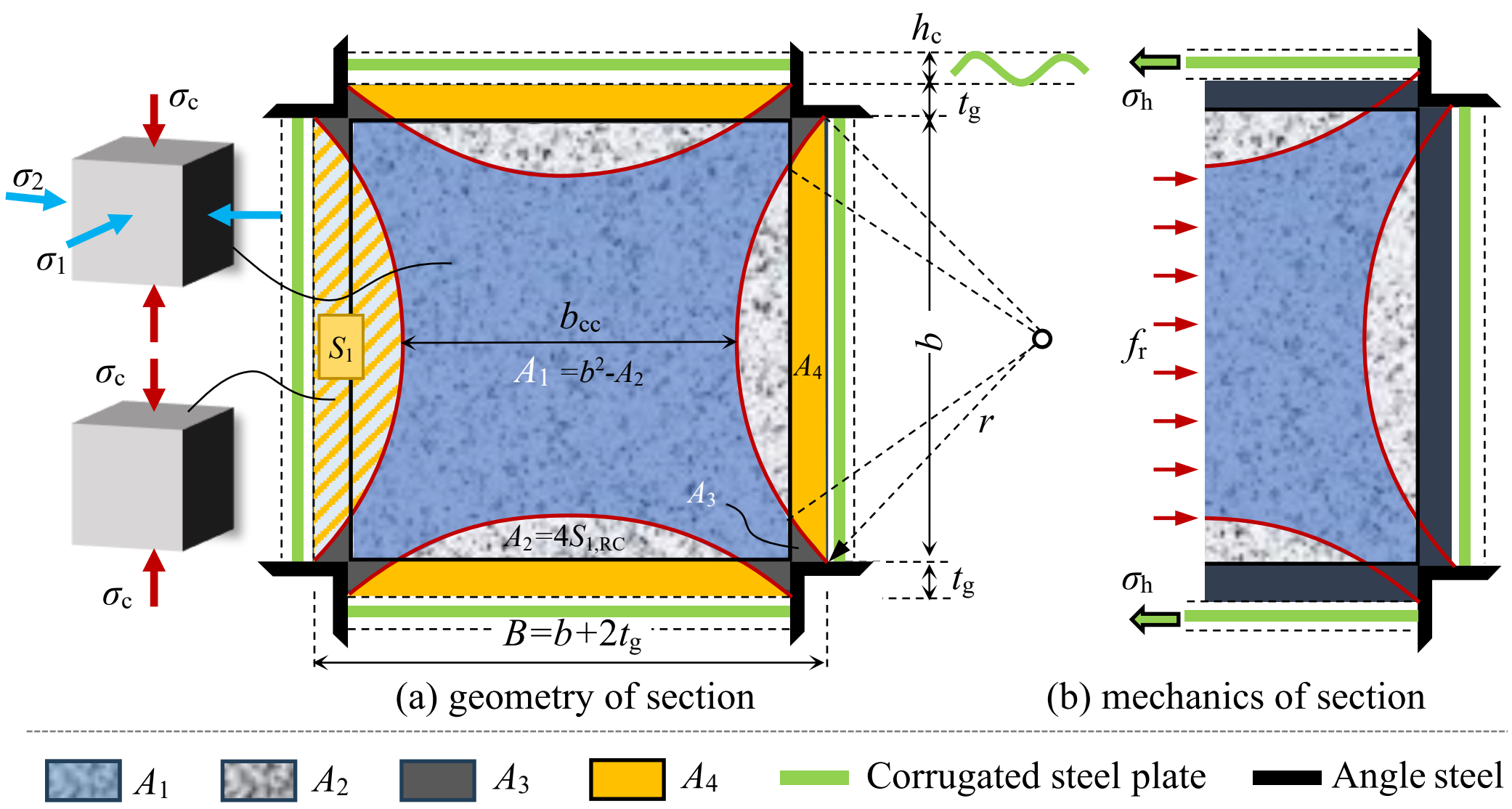

**Fig. 16. Internal force analysis of cross section.**

$$r = \frac{\sqrt{2}}{2} b \tag{15}$$

where the $B$, $b$, $r$, $h_c$ are width of enhanced RC section, width of grout layer, imaginary radius of unconfined zone, and corrugation depth.

**(1) Section geometry**

With the $r$, the single unconfined zone area and the unconfined RC section area can be easily calculated by **Eqs. (16)** and **(17)**.

$$S_1 = \frac{(\pi - 2) r^2}{4} \tag{16}$$

$$S_{1,RC} = r^2 \arccos\left(\frac{\sqrt{2}}{2} + \frac{t_g}{r}\right) - \left(\frac{\sqrt{2}}{2} r + t_g\right)\sqrt{\frac{r^2}{2} - \sqrt{2} t_g r - t_g^2} \tag{17}$$

As the partition shown in **Fig. 15**, the corresponding area can be calculated by using **Eqs. (18) – (21)**, where the edges between unconfined and confined area is assumed to be arc lines.

$$A_1 = b^2 - A_2 \tag{18}$$

$$A_2 = 4S_{1,RC} \tag{19}$$

$$A_4 = 4(S_1 - S_{1,RC}) \tag{20}$$

$$A_3 = 4(B^2 - t_g^2) - b^2 \tag{21}$$

**(2) Strength models**

For the A1 region, the concrete core is confined in a damaged state and thus its compressive strength should account for corresponding plastic damage. In this confined region, the confined concrete compressive strength should be taken into analysis of load bearing capacity.

There are plenty of strength model for confined concrete, which is primary related to the confining stress.

**Table 4. Strength models for confined concrete.**

| | | |
|---|---|---|
| Type I: $f_{cc}/f_{co} = a + k\left(f_r/f_{co}\right)^b$ | | |
| Park R. et al. [57] | square confined concrete | $f_{cc} = f_{co} + 2f_r$ |
| Hoang L. et al.[58] | circular steel tube confined concrete | $f_{cc} = f_{co} + 3.16f_r$ |
| Richart E.[59] | combined compressed concrete | $f_{cc} = f_{co} + 4.1f_r$ |
| Eurocode 2 [60] | concrete filled steel tube | $f_{cc} = f_{co} + 5f_r$, $f_r < 0.05f_{co}$ |
| Zhou X. et al. [61] | square confined concrete | $f_{cc} = f_{co} + 5.1f_r$ |
| Candappa D. et al.[62]<br>Wang Y. et al.[33] | concrete under triaxial compression | $f_{cc} = f_{co} + 5.3f_{co}$ |
| Saatcioglu M. et al.[63] | confined concrete | $f_{cc} = f_{co} + 6.7f_r^{0.83}$ |
| Hoshikuma J. et al.[64] | confined concrete | $f_{cc} = f_{co} + 7.6f_r$ |
| Fafitis A. et al.[65] | confined concrete | $f_{cc} = f_{co} + \left(1.15 + \dfrac{21}{f_{co}}\right)f_r$ |
| Megahed K. et al. (2024) [66]<br>Lim. J. and Ozbakkaloglu T. (2014) [67] | CFST stub columns | $f'_{cc} = f'_c + 5.2f_c'^{0.91}\left(\dfrac{f_r}{f'_c}\right)^{f_c'^{-0.06}}$ |
| Type II: $f_{cc} = kf_{co}$ | | |
| Mander J. et al.[68] | stirrups confined concrete | $f_{cc} = f_{co}\left(-1.254 + 2.254\sqrt{1 + 7.94\dfrac{f_r}{f_{co}}} - 2\dfrac{f_r}{f_{co}}\right)$ |
| Attard M. [69] | confined concrete | $f_{cc} = f_{co}\left(\frac{f_r}{f_t}+1\right)^{1.25 f_c^{-0.21}(1+0.062 f_r/f_{co})}$ |
| Setunge S. et al.[70] | confined high-strength concrete | $f_{cc} = f_{co}\left(1+13.07 f_r/f_{co}\right)^{0.63}$ |
| Cusson D. et al.[71] | confined high-strength concrete | $f_{cc} = f_{co}\left[1+2.1\left(f_r/f_{co}\right)^{0.7}\right]$ |
| Legeron F. et al.[72] | confined normal & high-strength concrete | $f_{cc} = f_{co}\left[1+2.4\left(f_r/f_{co}\right)^{0.7}\right]$ |
| Bousalem B. et al.[73] | rectangular confined concrete | $f_{cc} = f_{co}\left(1+0.8\frac{f_r}{\sqrt{f_{co}}}\right)$ |
| Girgin Z. et al. [74] | confined ultra-high-strength concrete | $f_{cc} = f_{co}\left[1+4.08\left(\frac{f_r}{f_{co}}\right)^{0.83}\right]$ |

Note: the $f_{cc}$, $f_{co}$ and $f_r$ are the confined concrete compressive strength, compressive strength of plain concrete and confining stress, respectively.

The listed models provide two typical calculation modes for the concrete under triaxial compression, including linear and non-linear equations. However, the models are proposed based on different test conditions, which may lead to discreteness in predicting the strength of concrete confined by square corrugated steel jacket.

As proposed by Park R. et al. [57], the strength model of confined square concrete is modified for the pre-damaged concrete confined by square corrugated steel jacket, i.e., $f'_{cc}$ in **Eq. (22)**.

$$f'_{cc} = f'_c + 2f_r \tag{22}$$

where $f'_c$ and $f_r$ are the residual compressive strength of concrete with damage and confining stress.

As shown in **Fig. 16**, the $f_r$ can be obtained by the equilibrium equation between the resultant forces of corrugated steel tube and confined concrete in A1 section:

$$f_r = \frac{\eta_{cs} t_{cs} \sigma_h}{B-(2-\sqrt{2})b} \tag{23}$$

where $\eta_{cs}$, $t_{cs}$, $\sigma_h$, $B$, and $b$ are respectively the amplification coefficient [32], wall thickness and average horizontal stress of corrugated steel, width of jacket and width of RC section.

For the A2 region, the compressive strength of concrete with damage ($f'_c$) should be employed, and the influence of damage can be expressed as **Eq. (24)**. A strength degradation factor of the concrete with damage, $k_d$, can be utilized to characterize the reduced strength caused by the damage.

$$f'_c = f_c\left(1 - k_d\sqrt{\frac{\varepsilon_0}{\varepsilon_c}}\right) \tag{24}$$

$$k_{\mathrm{d}}=\begin{cases}-0.02\left(\varepsilon_0/\varepsilon_{\mathrm{c}}\right)^2+0.12\varepsilon_0/\varepsilon_{\mathrm{c}}\ , & 0<\dfrac{\varepsilon_0}{\varepsilon_{\mathrm{c}}}\le 3\\ -0.006\left(\varepsilon_0/\varepsilon_{\mathrm{c}}\right)^2+0.03\left(\varepsilon_0/\varepsilon_{\mathrm{c}}\right)+0.144\ , & 3<\dfrac{\varepsilon_0}{\varepsilon_{\mathrm{c}}}\le 7\end{cases} \tag{25}$$

where $\varepsilon_0$ is the compressive strain at the start of unloading, and the compressive strain in the concrete at the peak stress $\varepsilon_c=0.7f_c^{0.31}$ (EC2 [60]). For the strengthened concrete column suffering damage, the $f_c'$ should be employed to calculate the ultimate load bearing capacity.

As per the study on concrete-filled corrugated steel tube under cyclic triaxial compression (Fang & Wang et. al. [75],[76]), $k_d$ can be obtained between 0.03 to 0.20 at the ultimate compressive strength. Based on the experimental stress-strain curves of the concrete under cyclic compression, the $k_d$ at different loading stage is accumulated and the regression analysis is conducted to capture its change. To better describe the rule of the loss of strength in different damage stages, the initial damage ratio $\varepsilon_0/\varepsilon_c$ is proposed to further match the $f_c'$ with the corresponding damage. As shown in **Fig. 17**, it grows from 0 in the ascending stage of the compression loading, and it continues to increase. The $k_d$ only decreases until the compressive strain exceeds 3 times of the peak strain of the plain concrete without confinement or lateral pressure.

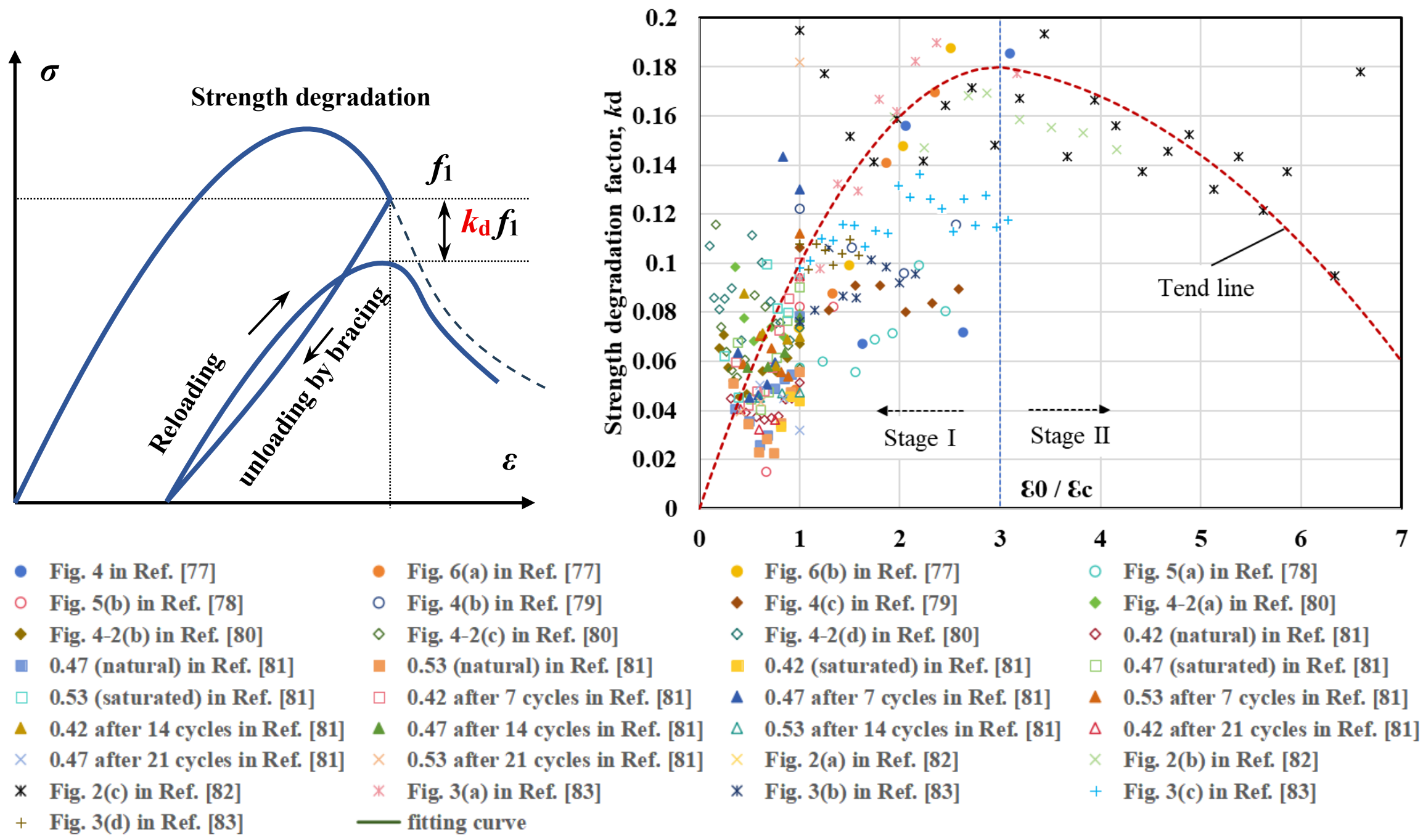

**Fig. 17. Strength degradation of concrete.**

Similarly, for the A3 zone, the infilled grout is confined without damaged before loading. Therefore, the confined compressive strength of grout ($f_{cc,g}$) should be adopted to calculate the load bearing capacity, which is given by **Eq. (26)**. For the A4 zone under less confinement, the compressive strength of grout ($f_{c,g}$) is directly employed to calculate the capacity.

$$f_{\mathrm{cc,g}} = f_{\mathrm{c,g}} + 2f_{\mathrm{r}} \tag{26}$$

where $f_{cc,g}$ and $f_{c,g}$ are confined reinforcement material compressive strength, and reinforcement material compressive strength, respectively.

## 5.2 Ultimate load-bearing capacity design methods

In this section, existing design specifications for concrete-filled steel tubular (CFST) and conventional concrete columns are adapted to establish a prediction model for the ultimate load-carrying capacity of axially compressed RC short columns strengthened with corrugated steel jackets. Considering the potential contribution of angle steels to the section stiffness at critical points, the steel properties of both the corrugated sheets and angle steels are incorporated into the overall steel ratio calculation of the equivalent steel tube. The yield strength of the equivalent steel tube is determined as the weighted average of the confining steel components, as presented in **Eq. (27)**.

$$f_b = \frac{f_{\mathrm{y,L}} A_L + f_{\mathrm{y,cs}} A_{\mathrm{cs}}}{A_{\mathrm{a}}} \tag{27}$$

where $f_{\mathrm{y,L}}$, $f_{\mathrm{y,cs}}$, $A_{\mathrm{a}}$, $A_L$, $A_{\mathrm{cs}}$ are the yield strength of angle steel and corrugated steel, cross sectional areas of equivalent steel tube ($A_{\mathrm{a}}=A_L+A_{\mathrm{cs}}$), and the area of angle steel and corrugated steel, respectively. The $f_{\mathrm{y,cs}}$ can be obtained by the strength of the peak and flat segment of corrugation (**Table. 2**) using weighted average value, with weighted coefficients of 0.3 and 0.7 [32], respectively.

**(1) Eurocode 4**

The BS EN 1994-1-1:2004 (Eurocode 4 [85]) provide a calculation method for axially loaded CFST columns, as given by **Eq. (28)**.

$$N_{\mathrm{u}} = \eta_{\mathrm{a}} f_{\mathrm{b}} A_{\mathrm{a}} + \bar{f}_{\mathrm{c}} A_{\mathrm{cc}} \left( 1 + \eta_{\mathrm{c}} \frac{t}{d} \frac{f_{\mathrm{b}}}{f_{\mathrm{ck}}} \right) + f_{\mathrm{y,s}} A_{\mathrm{s}} \tag{28}$$

where $f_{\mathrm{ck}}$ is the characteristic value of the cylinder compressive strength of concrete at 28 days, which was tested and calculated as MPa for the specimens in this study. $f_b$, $f_{y,\mathrm{s}}$, $A_{\mathrm{cc}}$, $A_{\mathrm{s}}$, $t$ and $d$ are the yield strength of confining steel and rebars, confined concrete core and rebars, thickness and diameter of steel tube, respectively. For the strengthened section in this study, the $f_{\mathrm{c}}$ in **Eq. (28)** should be obtained by weighted average according to the area, namely, $\overline{f}_{\mathrm{c}}=(f_{\mathrm{c}}A_{\mathrm{c}}+f_{\mathrm{c,g}}A_{\mathrm{g}})/(A_{\mathrm{c}}+A_{\mathrm{g}})$, where the $f_{\mathrm{c}}$, $f_{\mathrm{c,g}}$, $A_{\mathrm{c}}$ and $A_{\mathrm{g}}$ are the compressive strength of concrete and grout, and the areas of concrete and grout, respectively.

$$\eta_{\mathrm{a}}=0.25\left(3+2\overline{\lambda}\right) \tag{29}$$

$$\eta_{\mathrm{c}}=4.9-18.5\overline{\lambda}+17\overline{\lambda}^{2} \tag{30}$$

The relative slenderness is given by $\overline{\lambda}=\sqrt{N_{\mathrm{R}}/N_{\mathrm{cr}}}$. The characteristic value of the plastic resistance to compression with no consideration of confinement:

$$N_{\mathrm{R}}=f_{\mathrm{b}}A_{\mathrm{a}}+0.85f_{\mathrm{c}}A_{\mathrm{c}}+f_{\mathrm{y,s}}A_{\mathrm{s}} \tag{31}$$

For the elastic critical force ($N_{\mathrm{cr}}$), it can be calculated based on the effective flexural stiffness ($EI_{\mathrm{eff}}$):

$$N_{\mathrm{cr}}=EI_{\mathrm{eff}}\frac{\pi^{2}}{l^{2}} \tag{32}$$

$$EI_{\mathrm{eff}}=0.6E_{\mathrm{c}}I_{\mathrm{c}}+E_{\mathrm{s}}I_{\mathrm{s}}+E_{\mathrm{a}}I_{\mathrm{a}} \tag{33}$$

where $E_{\mathrm{a}}$, $E_{\mathrm{c}}$, $E_{\mathrm{s}}$, and $I_{\mathrm{a}}$, $I_{\mathrm{c}}$, $I_{\mathrm{s}}$, respectively, are the elastic modular and second moments of area of the structural steel section. When calculating the $E_{\mathrm{c}}I_{\mathrm{c}}$, the grout should be also considered for such a strengthened composite section.

**(2) ACI-318**

The United States' code ACI-318 [86] directly employs the **Eq. (31)** for the CFST section under axial compression, where the confinement is neglected.

**(3) AISC-360**

The ANSI/AISC 360-16 [87] use a similar equation as the ACI-318, while different coefficient (0.95) for the resistance of concrete section is adopted as following:

$$N_{\mathrm{u}}=f_{\mathrm{b}}A_{\mathrm{a}}+0.95f_{\mathrm{c}}A_{\mathrm{c}}+f_{\mathrm{y,s}}A_{\mathrm{s}} \tag{34}$$

**(4) AIJ-CFT 1997**

The Japanese code AIJ-CFT [88] also employs a formula consists of concrete and tube to calculate the ultimate load bearing capacity of square CFST column under axial compression. Considering the contribution of longitudinal bars, it can be given by **Eq. (35)**.

$$N_u = 1.27 f_b A_a + 0.85 f_c A_c + f_{y,s} A_s \quad (35)$$

**(5) GB 50936-2014**

In the Chinese code GB50936-2014 [89], the section is regarded as a composite section that an equivalent martial strength is defined to calculate the axial compressed capacity of CFST section, as expressed by **Eq. (36)**. Similarly, the capacity of rebars should be also taken into account. The calculation method for the confinement factor, $\xi$, should be modified for such a strengthened section, as its core consists of concrete and grout.

$$N_u = f_{sc} A_a + f_{y,s} A_s \quad (36)$$

$$f_{sc} = \left(1.212 + B_c \xi + C_c \xi^2\right) f_c \quad (37)$$

$$B_c = 0.723 + 0.615 \times 10^{-3} f_y \quad (38)$$

$$C_c = 0.026 - 4.861 \times 10^{-3} f_c \quad (39)$$

where $f_{sc}$ is the compressive strength of the composite section, and $\xi$ is the equivalent confinement factor of confining steel. $f_c$ and $f_{c,g}$ are the compressive strength of concrete and grout respectively. $A_c$ and $A_g$ are the areas of concrete and grout correspondingly.

**(6) GB/T51446-2021**

The expression in the **GB/T51446-2021** [90] is the same with **Eq. (36)**, where the nominal strength of the composite CFST section is replaced by **Eq. (40)**.

$$f_{sc} = \left(1.14 + 1.02\xi\right) f_c \quad (40)$$

**(7) BS 5400-5:2005**

For a square section, the British code BS5400-5:2005 [91] directly use the unconfined concrete strength, and the contribution of steel tube is taken into account. which can be expressed by **Eq. (41)** for the strengthened RC short column in this study.

$$N_u = 0.45\left(f_c A_c + f_{c,g} A_g\right) + 0.95\left(f_{y,L} A_L + f_{y,cs} A_{cs}\right) + 0.87 f_{y,s} A_s \quad (41)$$

where the subscripts 'c' and 'g' denote concrete and grout correspondingly, and the subscripts 'ac' and 'cs' represent angle steel and corrugated steel, respectively. The 's' means longitudinal bars.

**(8) GB 50010-2010**

By ignoring the load bearing capacity of corrugated steel, and neglecting its confinement effect,

the calculation method for a conventional RC cross section [92] can be also utilized to conservatively estimate the ultimate axial capacity of strengthened RC short column, see **Eq. (42)**.

$$N_u = 0.9\left(f_c A_c + f_{c,g} A_g + f_{y,L} A_L + f_{y,cs} A_{cs} + f_{y,s} A_s\right) \tag{42}$$

As compared in **Fig. 18**, the modified methods of design codes ACI-318, AISC-360, AIJ-CFT, BS5400 and GB50010-2010 can conservatively predict the ultimate load-bearing capacity of the axially compressed RC section strengthened with corrugated steel jacketing. The EC4, GB50936 and GB/T51446 provide overestimated capacity results, though.

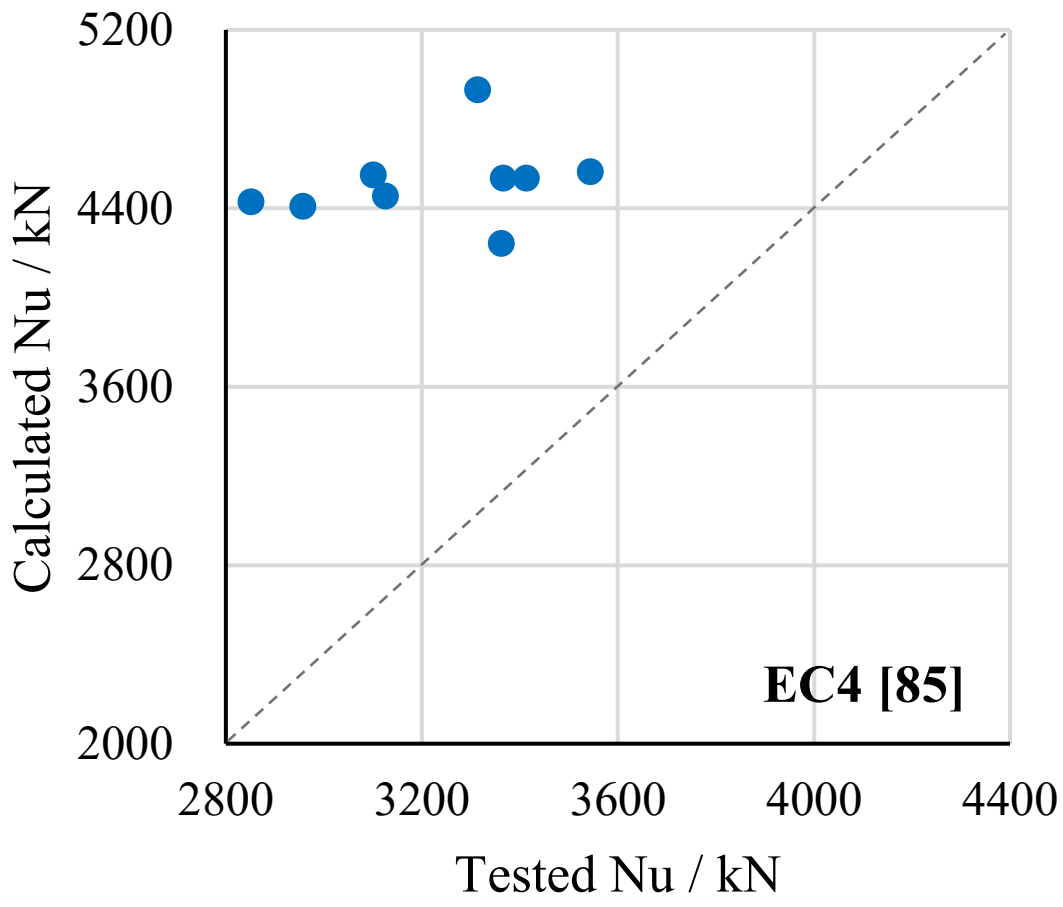

(a) EC4

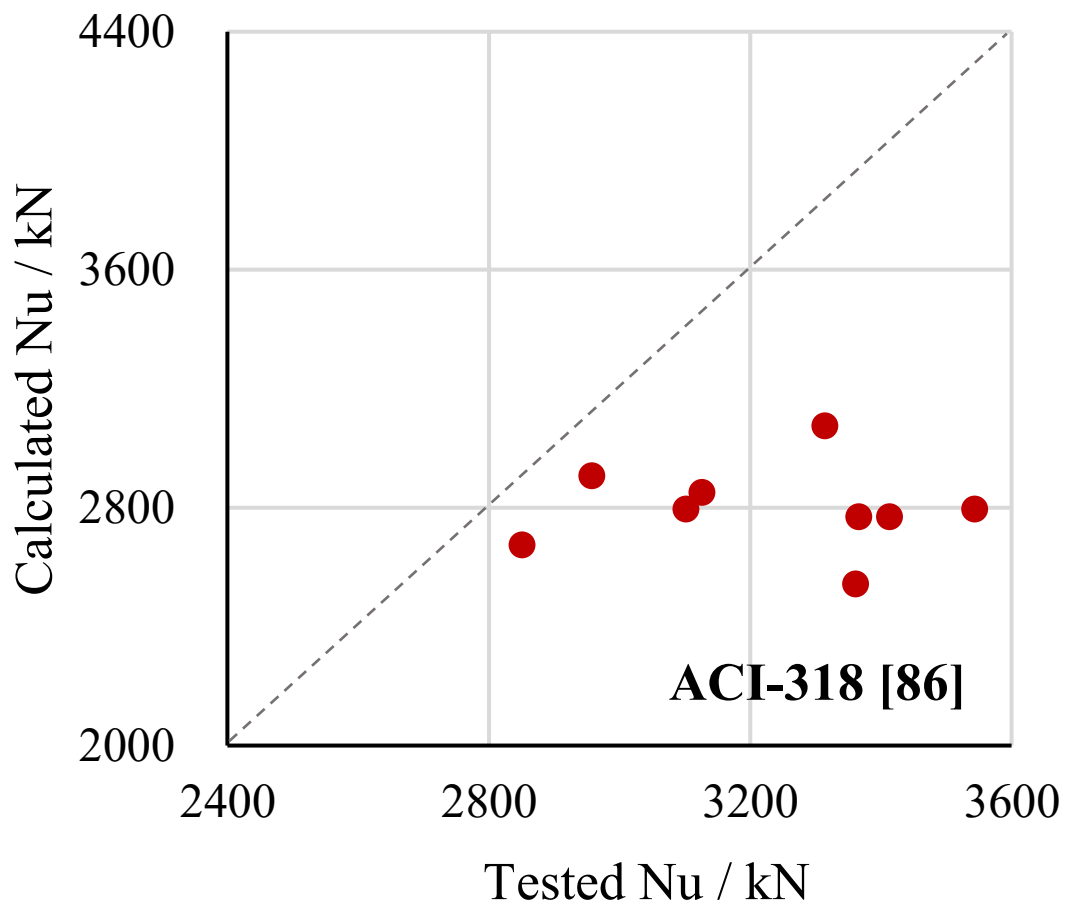

(b) ACI

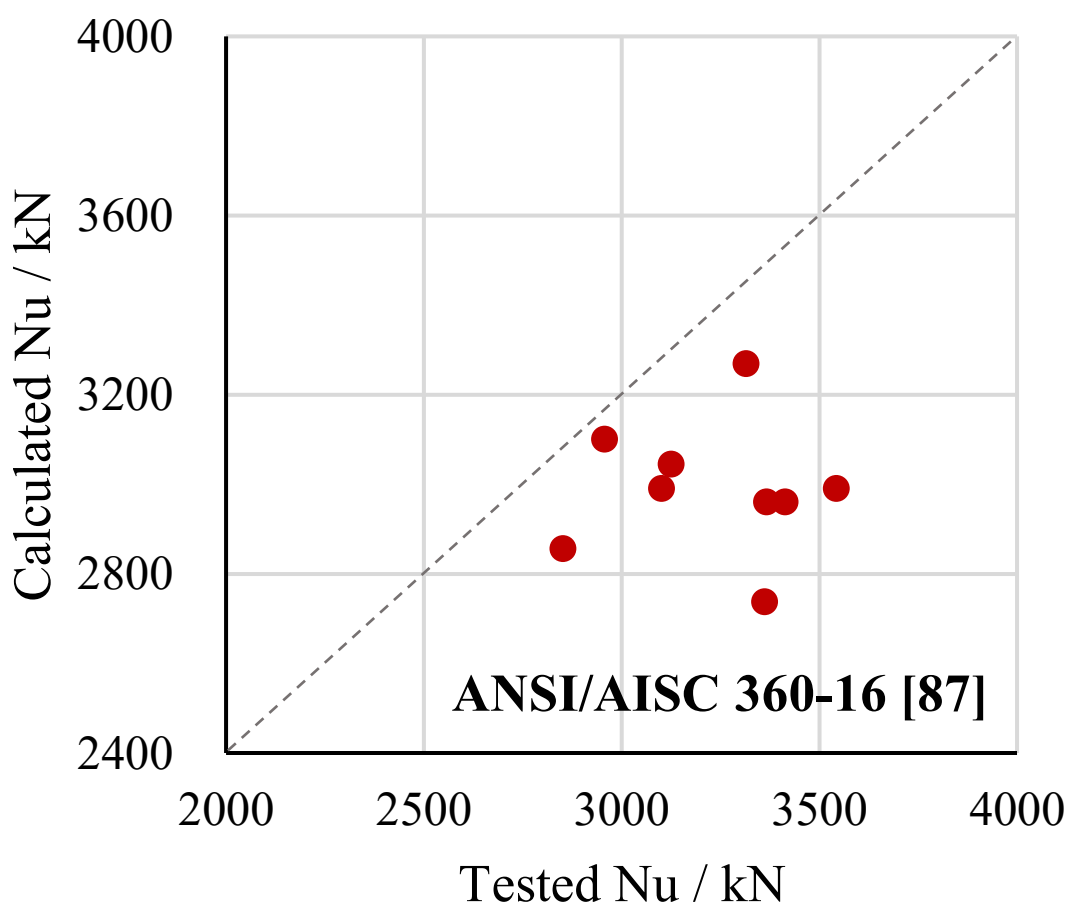

(c) AISC-360

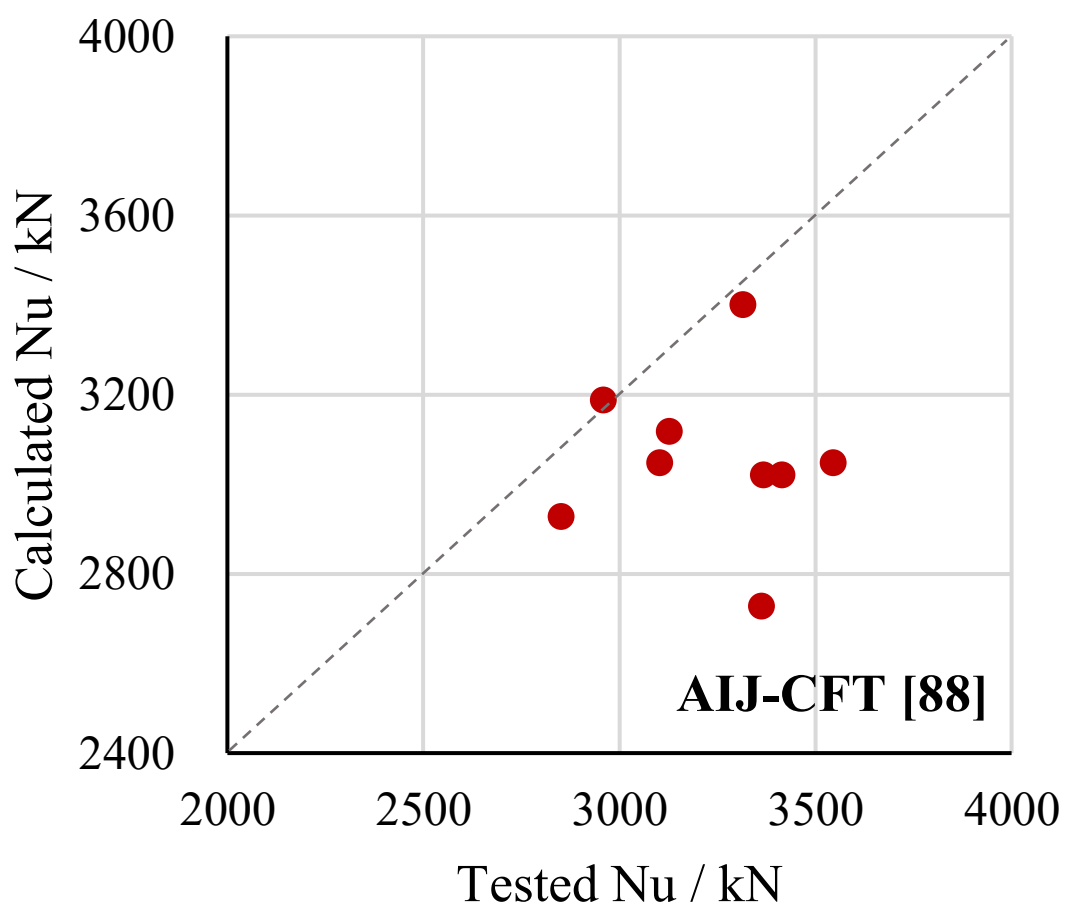

(d) AIJ-CFT

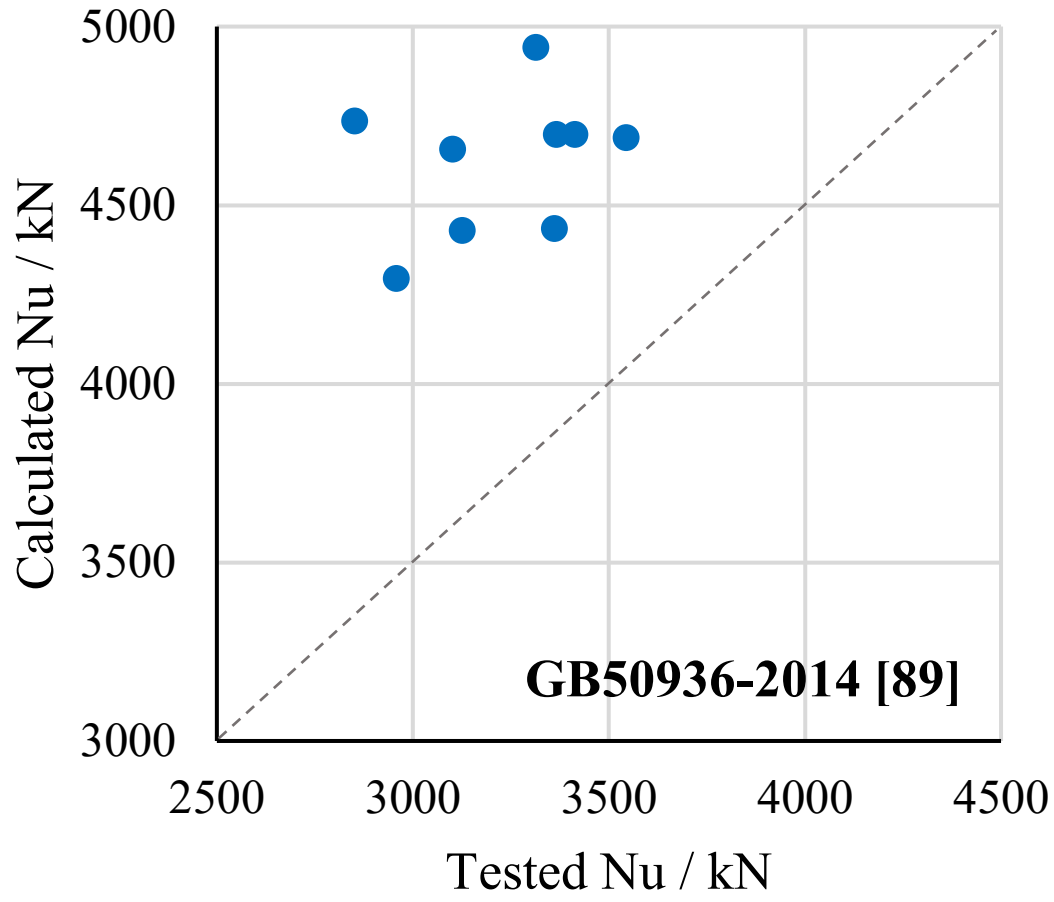

(e) GB50936

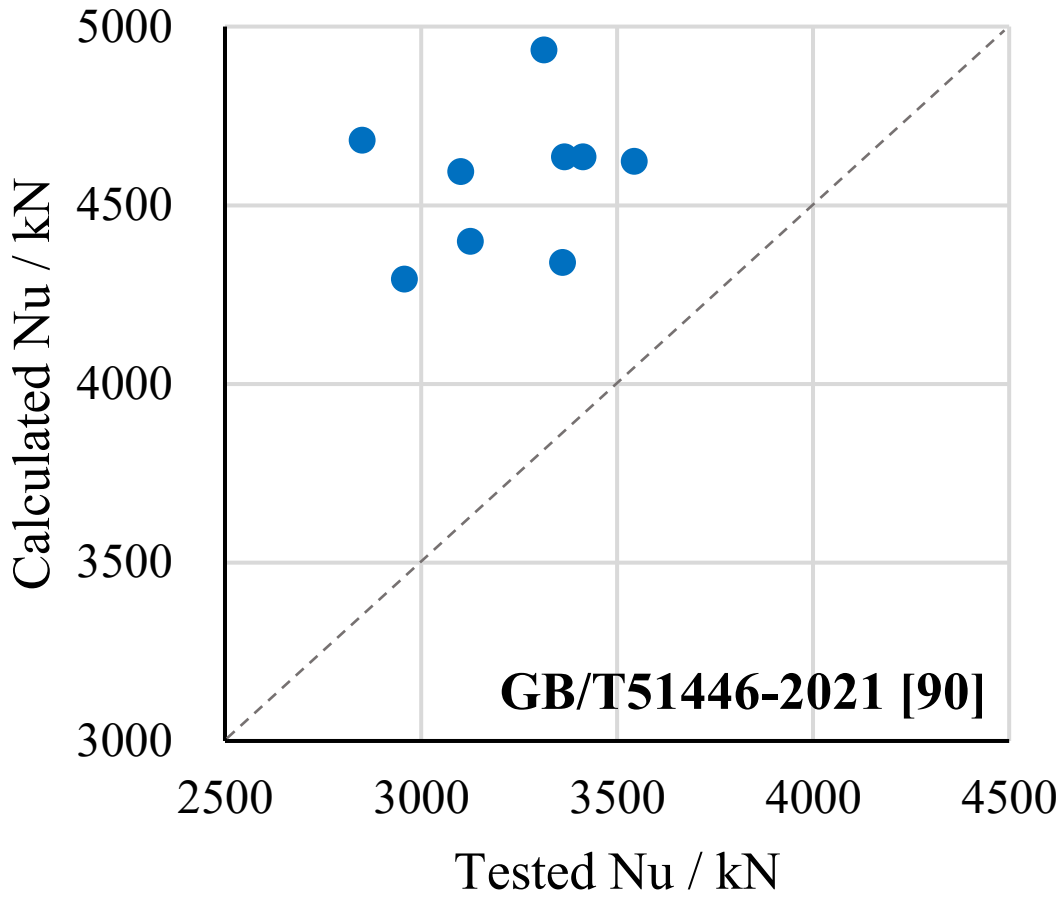

(f) GB/T51446

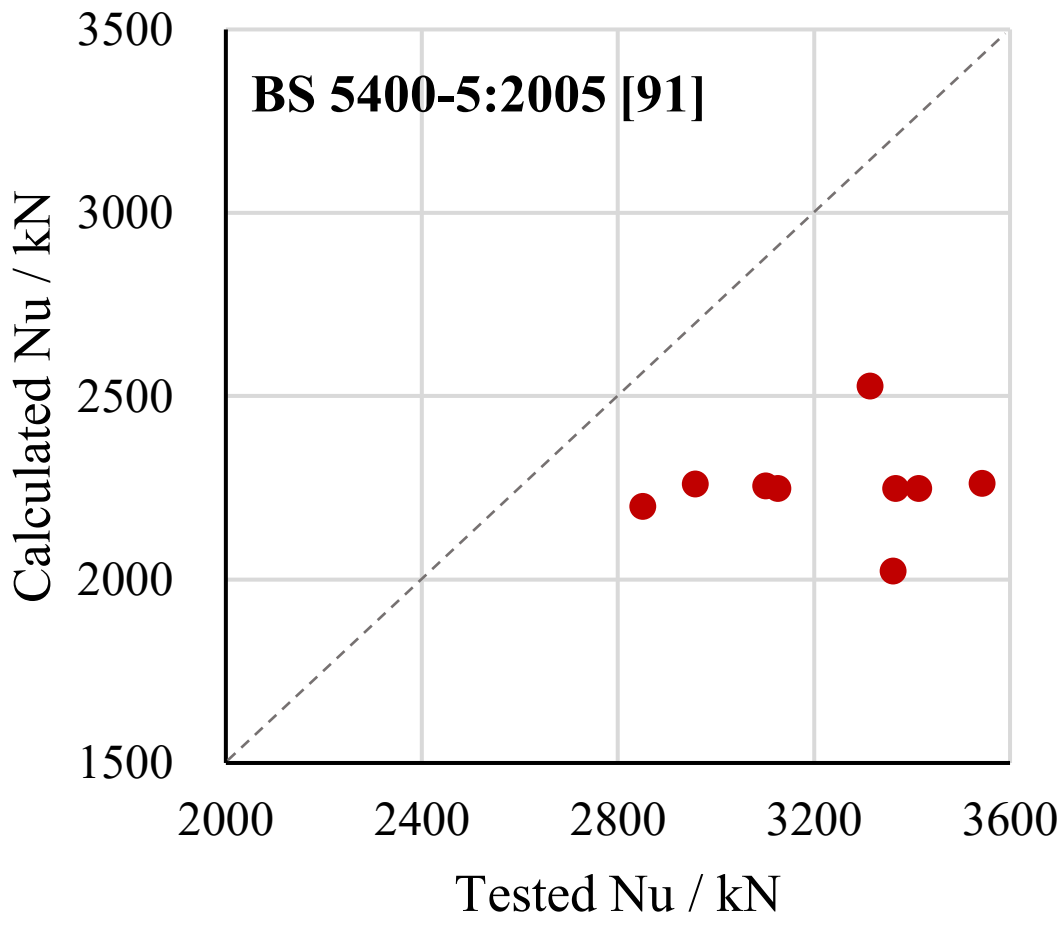

(g) BS5400

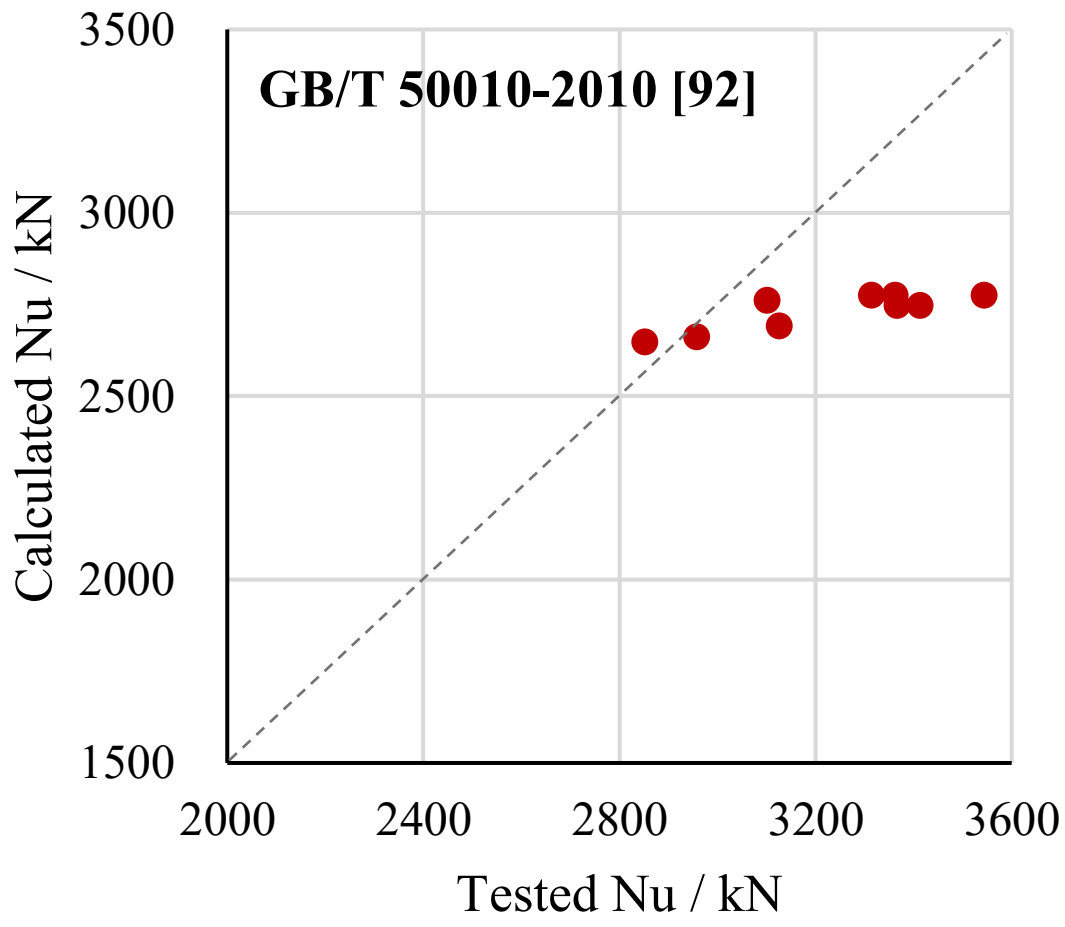

(h) GB/T 50010

**Fig. 18. Comparison of design codes.**

**(9) Proposed model**

Actually, the corrugated steel jacketed square RC column can be regarded as a type of confined concrete column [84]. As proven above, the corrugated steel sheet, especially near the centerline of cross section, does not directly bear axial load. Therefore, only the corner angle steels, reinforcing layer grout, core concrete, and reinforcement bars are considered to carry the load. The axial compressive load-bearing capacity of the strengthened RC section ($N_u$) can be calculated by:

$$N_{\mathrm{u}} = f'_{\mathrm{cc}}A_1 + f'_{\mathrm{c}}A_2 + f_{\mathrm{cc,g}}A_3 + f_{\mathrm{c,g}}A_4 + f_{\mathrm{y,L}}A_{\mathrm{L}} + f_{\mathrm{yb}}A_{\mathrm{b}} \tag{43}$$

where the $f_{\mathrm{y,L}}$, $f_{\mathrm{yb}}$, $A_{\mathrm{L}}$, $A_{\mathrm{b}}$ are yield strength of angle steel and original longitudinal rebars, areas of angle steel and longitudinal rebars, respectively.

Due to the RC section encapsulated with corrugated sheets possessing characteristics of confined concrete, the confined concrete strength model is applicable to the RC column rehabilitated with corrugated jacketing. However, if without consideration of the strength degradation, only substituting $f_{cc}$ into **Eq. (43)** cannot accurately predict the ultimate bearing capacity (**Fig. 19**). Therefore, the strength of confined pre-damaged concrete, $f'_{cc}$, should be adopted.

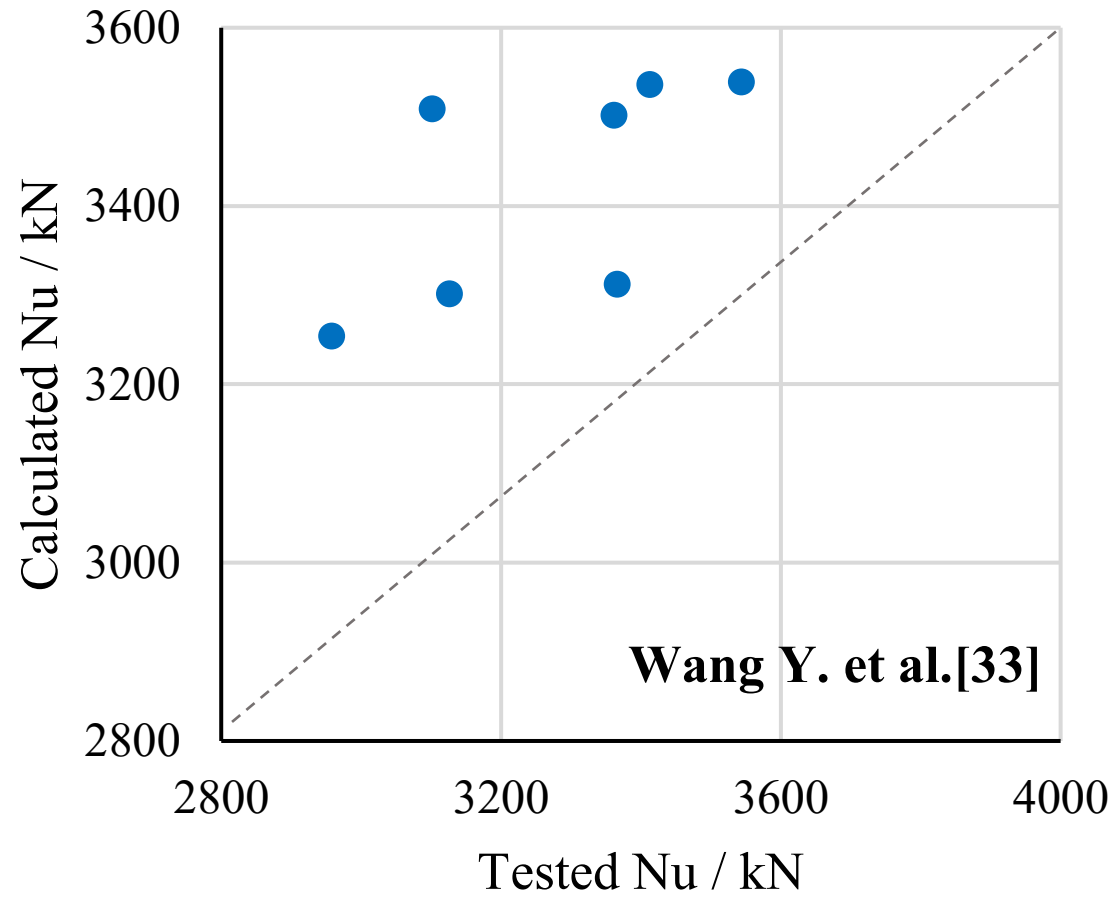

(a) concrete-filled corrugated steel tube model

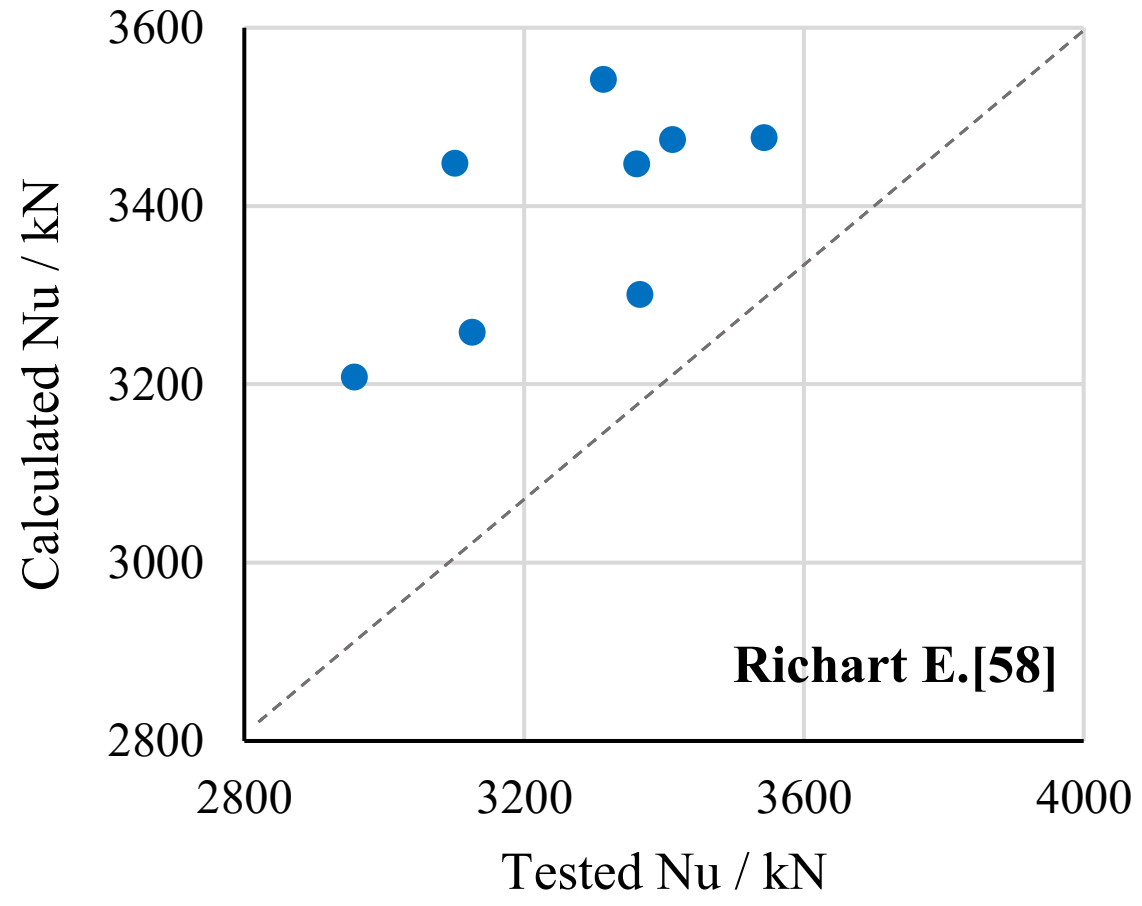

(b) hydraulic pressure confined concrete model

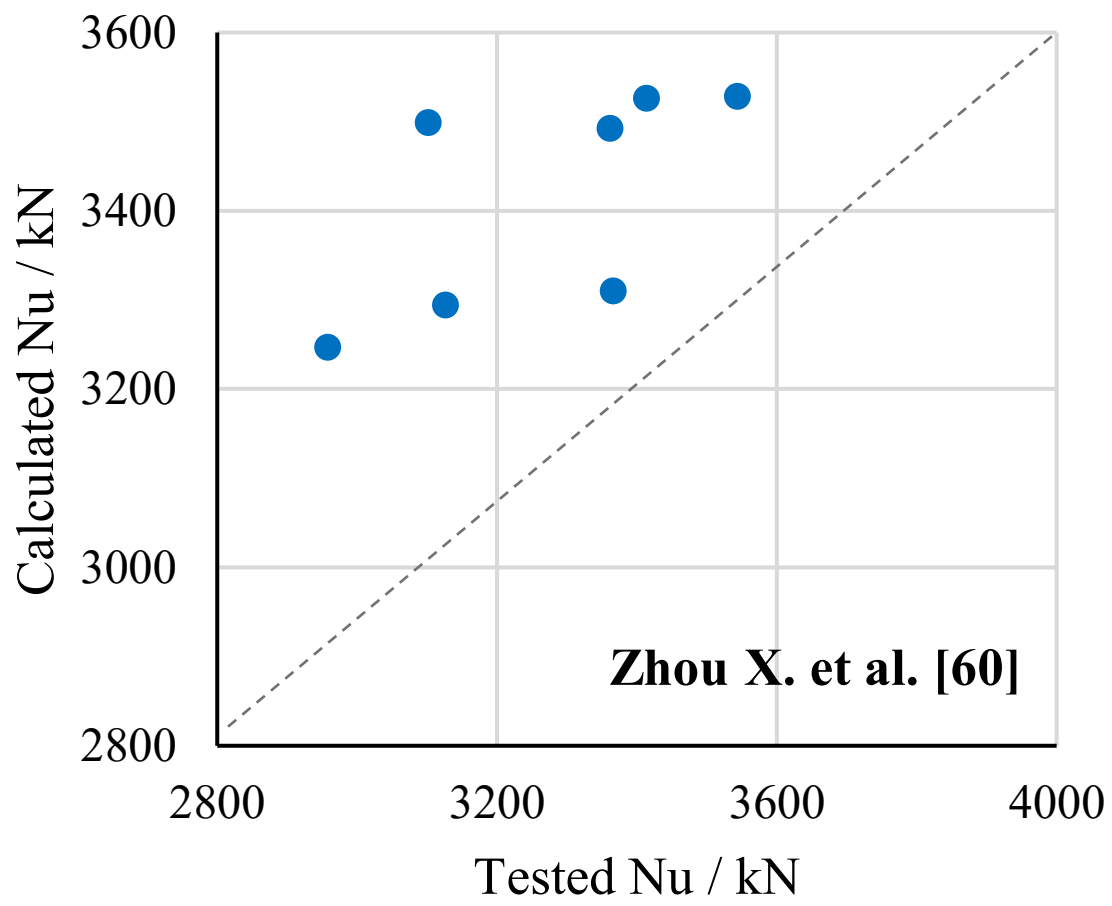

(c) square CFST strength model

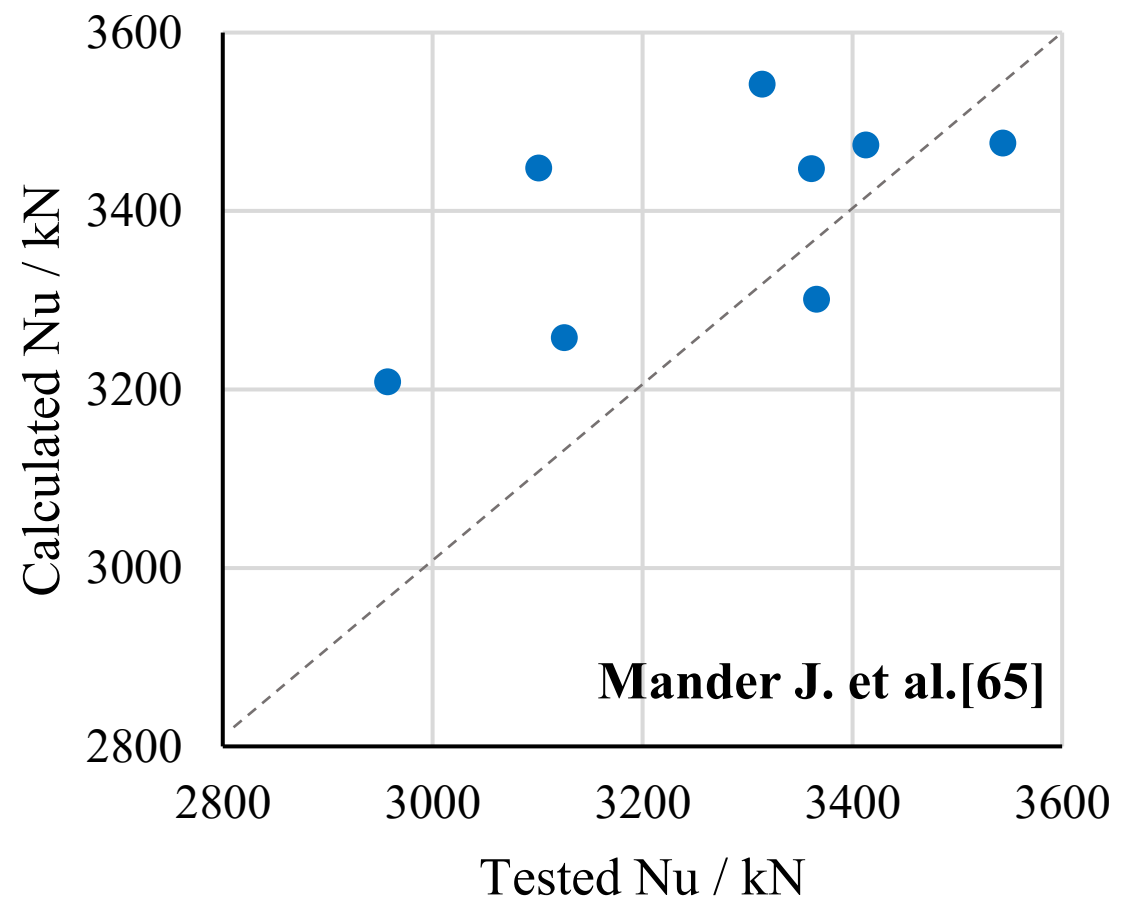

(d) spiral stirrups confined concrete model

**Fig. 19. Different strength models for load bearing capacity, $N_u$.**

As validated in **Fig. 20**, considering the contributions of each component to bearing capacity, the damage in the concrete core, and the enhancing effect of confinement on the concrete strength, the results calculated by using **Eq. (43)** closely matches the test piece actual measurements, with an average error of 0.4 %. It provides an accurate calculation method for predicting the ultimate bearing capacity of the reinforced section.

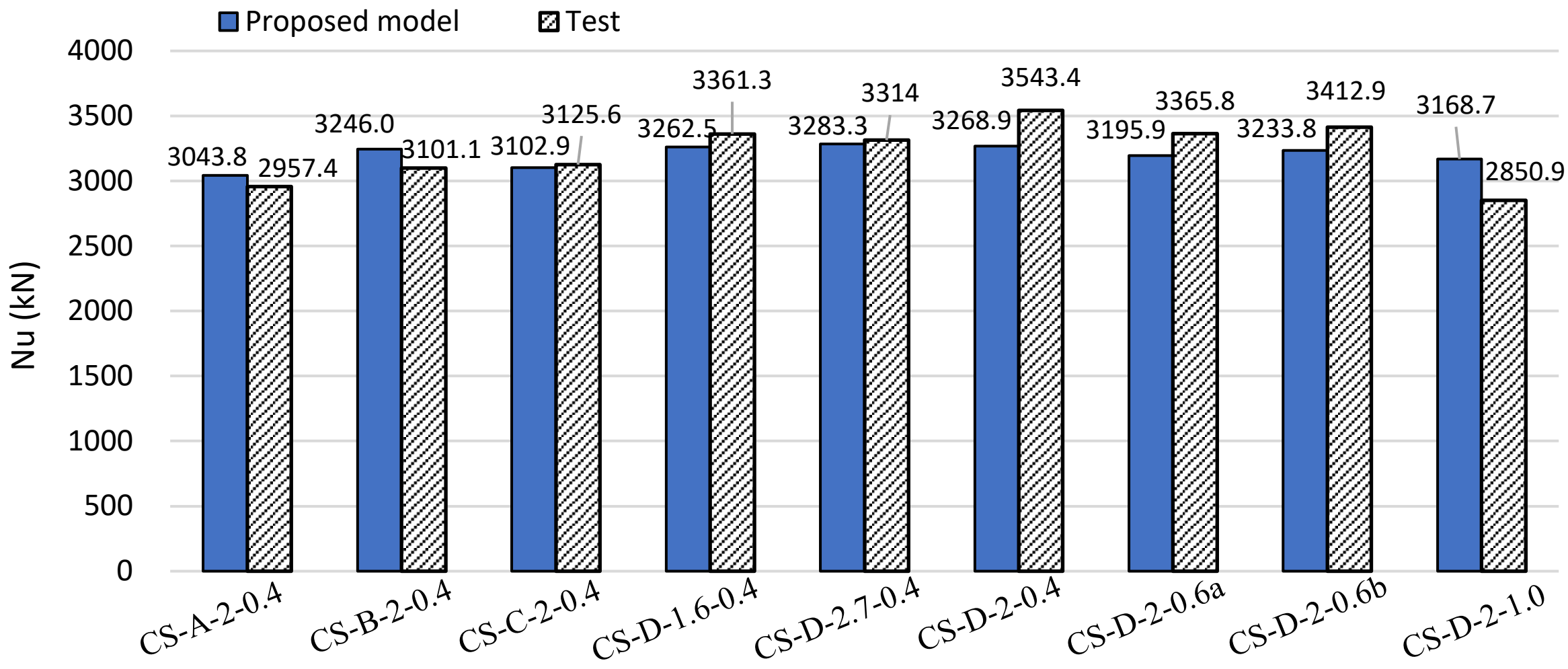

**Fig. 20. Validation of ultimate capacity model in this study.**

## 6. Conclusions

In this study, the axial behavior of pre-damaged RC short columns strengthened with corrugated steel jackets was experimentally investigated. Ten specimens were tested to examine the effects of connection methods, jacket thickness, and damage level on load–displacement response, stress distribution, and strain development. A strength model for confined concrete was developed, incorporating the influence of axial preload and transverse confinement. The main conclusions are as follows:

(1) Strengthening square RC columns with corrugated steel jackets significantly enhances structural performance. The ultimate bearing capacity of the jacketed reinforced concrete columns increased from 2117.7 kN to between 2850.9 kN and 3543.4 kN, corresponding to an improvement of 34.6% to 67.3% compared to the reinforced concrete section without jacketing.

(2) Due to limited confinement effectiveness in square sections, variations in jacket thickness had minimal influence on load–displacement behavior. The notable gains in ultimate capacity and ductility suggest that thinner jackets may be sufficient for strengthening square columns, provided construction stability is maintained.

(3) Increased preload ratios ($N_0/N_{u,RC}$) reduced the residual strength of the concrete, thereby diminishing the load-bearing capacity of the strengthened specimens. Contrary to conventional expectation, however, specimens subjected to higher $N_0/N_{u,RC}$ levels exhibited improved ductility.

From a mechanistic perspective, in a jacketed column after unloading, a higher $N_0/N_{u,RC}$ ratio induces greater concrete softening and enhanced lateral expansion under re-compression. This dilation promotes tighter interaction between the core concrete and the steel jacket, activating stronger passive confinement. Furthermore, in such a composite system, the reduction in concrete strength due to preloading conversely increases the confinement factor ($\xi$), which contributes to the improved ductility performance.

(4) The connection type of the jacket had minimal influence on the ultimate capacity or ductility of the strengthened columns. However, the specimen with fully welded corners provided the strongest confinement, achieving an ultimate bearing capacity of 3361.3 kN, 9.8% higher than the average of the other three specimens. Nevertheless, bolted connections remain a viable option in strengthening engineering due to their ease of assembly.

(5) In square sections, the confining potential of the corrugated steel jacket was not fully activated until concrete damage reached an advanced stage. In the severely damaged specimen (CS-D-2-1.0), the horizontal strain approached the yield point, and the corresponding horizontal stress reached approximately 200 MPa, far greater than those specimens with less damage. This indicates that greater concrete pre-damage elicits stronger passive confinement from the jacket, reflected in a higher confinement factor. Meanwhile, the longitudinal stresses in the jacket remained low, suggesting that its primary role is to provide lateral confinement rather than to resist axial loads directly.

(6) The existing confined concrete models are inadequate for accurately predicting the ultimate capacity of jacketed columns, as they typically overlook the influence of confinement effects and, in particular, fail to account for strength degradation resulting from pre-existing damage. It is therefore recommended to adopt a calculation method that incorporates these two critical factors to ensure conservative and reliable design of corrugated steel-jacketed RC columns.

## Acknowledgements

The study is supported by the Natural Science Foundation of China (No. 52308303), the Youth Foundation of Chongqing Construction Engineering Group, and the Research and Innovation Program for Graduate Students in Chongqing (Grant Nos. CYS25491). The financial supports are gratefully acknowledged.